\documentclass[preprint,11pt,3p,times]{article}

\usepackage{subcaption}
\usepackage{amsmath}
\usepackage{amssymb}
\usepackage{url}
\usepackage{authblk} 
\usepackage[left=2.2cm,right=2.2cm,top=2.0cm,bottom=2.0cm]{geometry}

\usepackage{tikz}
\usetikzlibrary{shapes,arrows, shapes.geometric,calc, snakes, spy, arrows.meta}

\tikzstyle{decision} = [diamond, draw, fill=gray!20, 
text width=11em, text badly centered, node distance=4cm, inner sep=0pt]
\tikzstyle{smalldecision} = [diamond, draw, fill=gray!20, 
text width=7em, text badly centered, node distance=3.5cm, inner sep=0pt]
\tikzstyle{block} = [rectangle, draw, fill=gray!20, 
text width=5em, text centered, rounded corners, minimum height=4em, node distance=4cm]
\tikzstyle{longblock} = [rectangle, draw, fill=gray!20, 
text width=21em, text centered, rounded corners, minimum height=4em, node distance=4cm]
\tikzstyle{line} = [draw, -latex']
\tikzstyle{cloud} = [draw, ellipse,fill=gray!20, text width= 9em, text centered, node distance=7cm,
minimum height=4em]

\pgfmathsetmacro{\radius}{0.25}
\pgfmathsetmacro{\length}{5}

\tikzset {
	bent cylinder/.pic={
		\pgfmathsetmacro{\anglebottom}{90-20+40*rnd}
		\pgfmathsetmacro{\strengthbottom}{\length*(0.3+rnd/3)}
		\pgfmathsetmacro{\angletop}{90-20+40*rnd}
		\pgfmathsetmacro{\strengthtop}{\length*(0.3+rnd/3)}
		
		\fill[fill=white] 
		(\radius,0) 
		.. controls +(\anglebottom:\strengthbottom) 
		and +(-\angletop:\strengthtop) 
		.. (\radius,\length) 
		arc (360:180:\radius cm and 0.5*\radius cm) 
		.. controls +(-\angletop:\strengthtop) 
		and +(\anglebottom:\strengthbottom) 
		.. (-\radius,0) 
		;
		
		\draw[gray!90, ultra thick, dashed, line join=round]
		(\radius,0) 
		.. controls +(\anglebottom:\strengthbottom) 
		and +(-\angletop:\strengthtop) 
		.. (\radius,\length) 
		arc (360:180:\radius cm and 0.5*\radius cm) 
		.. controls +(-\angletop:\strengthtop) 
		and +(\anglebottom:\strengthbottom) 
		.. (-\radius,0) 
		;
		
		\draw[gray!90, ultra thick, line join=round]
		(0,0) 
		.. controls +(\anglebottom:\strengthbottom) 
		and +(-\angletop:\strengthtop) 
		.. (0,\length) 
		;
		
		\filldraw[ultra thick, dashed, draw=gray!90, fill=white]
		(0,\length) circle (\radius cm and 0.5*\radius cm);
		
		\draw[ultra thick, dashed, draw=gray!90]
		(0,0) circle (\radius cm and 0.5*\radius cm);
	}
}

\definecolor{colUniBwOr}{rgb}{0.929,0.431,0.0} 

\graphicspath{{pictures/twoway_oscillating/}{pictures/biegebalken_comparison/}{pictures/plantpatch/very_new/}{pictures/solver_comparison/biegebalken/}{pictures/biegebalken_comparison/new/}{pictures/solver_comparison/biegebalken_light/}}

\newtheorem{problem}{Problem}
\newtheorem{remark}{Remark}

\newcommand{\indFBI}{\mathrm{fbi}}
\newcommand{\indFluid}{\mathrm{f}}
\newcommand{\indBeam}{\mathrm{b}}
\newcommand{\indLagMult}{\mathrm{lm}}
\newcommand{\trans}{\mathrm{T}}
\newcommand{\ndim}{n^{\mathrm{dim}}}

\newcommand{\indLagMultDof}{p}
\newcommand{\indBeamDof}{q}
\newcommand{\indFluidDof}{r}

\newcommand{\fvelvec}{\mathbf{\hat{v}}_h^\indFluid}
\newcommand{\fvelnod}{\mathbf{\hat{v}}_h^{\indFluid,i}}
\newcommand{\fvelfunch}{\mathbf{v}_h^\indFluid}

\newcommand{\fpresvec}{\mathbf{\hat{p}}_h^\indFluid}
\newcommand{\fpresnod}{\hat{p}_h^{\indFluid, i}}
\newcommand{\fpresfunch}{p_h^\indFluid}

\newcommand{\rvec}{\mathbf{\hat{r}}_h}
\newcommand{\rveck}{\mathbf{\hat{r}}_h^k}
\newcommand{\rnodtrans}{\hat{\mathbf{d}}^j_h}
\newcommand{\rnodtang}{\hat{\mathbf{t}}^j_h}
\newcommand{\rfunch}{\mathbf{r}_h}

\newcommand{\bvelvec}{\mathbf{\hat{v}}_h^\indBeam}
\newcommand{\bvelveck}{\mathbf{\hat{v}}_h^{\indBeam, k}}

\newcommand{\bvelfunch}{\mathbf{v}_h^\indBeam}

\newcommand{\lambdavec}{\hat{\boldsymbol{\mathsf{\lambda}}}_h}
\newcommand{\lambdanod}{\hat{\boldsymbol{\mathsf{\lambda}}}_h^k}
\newcommand{\lambdafunch}{\boldsymbol{\mathsf{\lambda}}_h}

\newcommand{\brhsvec}{\mathbf{f}_h^\mathcal{B}}
\newcommand{\frhsvec}{\mathbf{f}_h^\mathcal{F}}
\newcommand{\btrhsvecn}{\mathbf{f}_{h}^{\mathcal{B}, n}}

\newcommand{\bfbiforcevec}{\mathbf{f}_h^{\mathcal{B}, \indFBI}}
\newcommand{\bfbiforceveck}{\mathbf{f}_h^{\mathcal{B}, \indFBI, k}}
\newcommand{\bfbiforcevectemplate}[1]{\mathbf{f}_h^{\mathcal{B}, \indFBI, #1}}
\newcommand{\bfbiforceveckold}{\mathbf{f}_h^{\mathcal{B}, \indFBI, k-1}}
\newcommand{\bfbiforcevecnewton}[2]{\mathbf{f}_h^{\mathcal{B}, \indFBI, #1, #2}}

\newcommand{\fbiforceveck}{\mathbf{f}_h^{\indFBI, k}}

\newcommand{\fbirhsk}{\mathbf{r}^{\indFBI,k}}
\newcommand{\fbirhstemplate}[1]{\mathbf{r}^{\indFBI,#1}}
\newcommand{\fbirhsnewton}[2]{\mathbf{r}^{\indFBI,#1, #2}}
\newcommand{\jac}{\tilde{\mathbf{J}}_{\indFBI}^{k}}

\newcommand{\ie}{i.e.}

\newcommand{\wrt}{w.r.t}

\begin{document}

\title{A fully coupled regularized mortar-type finite element approach for embedding one-dimensional fibers into three-dimensional fluid flow}

\author[1]{Nora Hagmeyer}

\author[1,2]{Matthias Mayr*}

\author[1]{Alexander Popp}

\affil[1]{\small{Institute for Mathematics and Computer-Based Simulation, University of the Bundeswehr Munich, Germany}}

\affil[2]{\small{Data Science \& Computing Lab, University of the Bundeswehr Munich, Germany}}


\date{}

\maketitle

\abstract[Summary]{The present article proposes a partitioned Dirichlet-Neumann algorithm,
that allows to address unique challenges arising from a novel mixed-dimensional coupling of very slender fibers embedded in fluid flow
using a regularized mortar-type finite element discretization.
The fibers are modeled via one-dimensional (1D) partial differential equations based on geometrically exact nonlinear beam theory,
while the flow is described by the three-dimensional (3D) incompressible Navier-Stokes equations.
The arising truly mixed-dimensional 1D-3D coupling scheme constitutes a novel approximate model and numerical strategy,
that naturally necessitates specifically tailored solution schemes to ensure an accurate and efficient computational treatment.
In particular, we present a strongly coupled partitioned solution algorithm based on a Quasi-Newton method for applications
involving fibers with high slenderness ratios that usually present a challenge with regard to the well-known added mass effect.
The influence of all employed algorithmic and numerical parameters, namely the applied acceleration technique,
the employed constraint regularization parameter as well as shape functions, on efficiency and results of the solution procedure is studied through appropriate examples.
Finally, the convergence of the two-way coupled mixed-dimensional problem solution under uniform mesh refinement is demonstrated,
a comparison to a 3D reference solution is performed,
and the method's capabilities in capturing flow phenomena at large geometric scale separation is illustrated by the example of a submersed vegetation canopy.

}

\section{Introduction}
\label{sec:intro}

The application of structural models with reduced spatial dimensionality to adequately describe complex mechanical behavior of slender bodies has a long-standing history. Even though immersed fiber-like structures appear in a multitude of applications, the use of such highly efficient models, as they arise from geometrically exact beam theory presented in~\cite{simo1985,Reissner1972, Meier2014, Meier2015, meier_contact} for example, is still not well-studied in the context of multi-physics problems.
One such application, in which immersed fibers have a significant impact on the fluid flow, is the modeling of optimal flow control. In~\cite{kunze2012}, an experimental analysis of different coatings of hairy flaps to control vertex shedding behind an immersed cylinder are reported, and in~\cite{favier2009} a homogenized model is used to model the hairy coating.
A similarly global effect on the flow can be observed in the case of terrestrial canopies and submerged vegetation. Even though one immersed fiber does not affect the fluid flow considerably, the movement and interaction of a large patch of fibers has a significant impact on the vorticity of the overall flow~\cite{Tschisgale2020, tschisgale2021, oconnor2019, oconnor2022}. O'Connor et al.~\cite{oconnor2019, oconnor2022} have recently studied this effect via surface-coupling of two-dimensional structures to the two-dimensional Navier-Stokes equations.
The employed fully resolved approach comes with the difficulty that a high mesh resolution is required to be able to resolve the slender fibers. This leads to a computationally expensive model, which is problematic to extend to three dimensions, making it a prime target for the employment of one-dimensional beam theory and the application of a mixed-dimensional model, as done in~\cite{wang2019, tschisgale2021}. Such mixed-dimensional models are more accurate than homogenized methods, and at the same time computationally more efficient than fully resolved models. Still, depending on the employed coupling approach, also mixed-dimensional models can lead to very different levels of accuracy and complexity in terms of degrees of freedom.

\begin{figure}
	\begin{subfigure}{0.3\textwidth}
		\includegraphics[width=\textwidth]{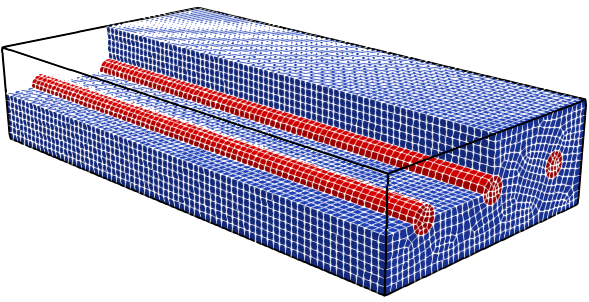}
	\end{subfigure}
	\begin{subfigure}{0.3\textwidth}
		\includegraphics[width=\textwidth]{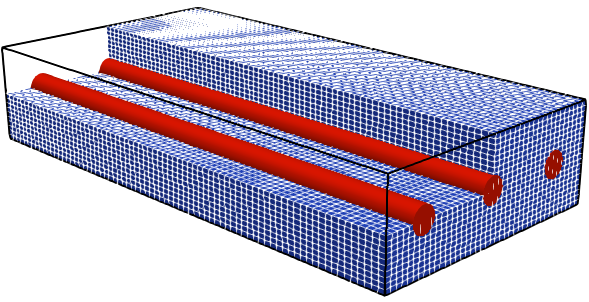}
	\end{subfigure}
	\begin{subfigure}{0.3\textwidth}
		\includegraphics[width=\textwidth]{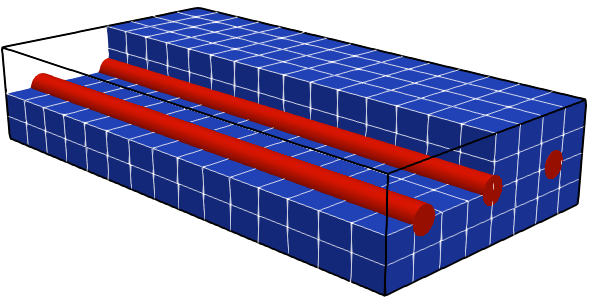}
	\end{subfigure}
	\caption{Various coupling strategies}
	\label{fig:coupling}
\end{figure}

Figure~\ref{fig:coupling} shows various modeling approaches for fluid-beam coupling along with the typically required resolution of the fluid mesh for each approach.
In all cases, the fibers are assumed to be modeled via a one-dimensional model.
The image on the left shows the reconstruction of the beams' actual surface and a subsequent fictitious meshing similar to the coupling used in~\cite{huang2019, Henshaw2015, ausas2022}. The fluid mesh resolves the fiber surfaces explicitly, either with matching or non-matching grids.
The middle picture shows a mesh with approximately 3-5 fluid elements over the beam diameter as commonly required for immersed type methods, for which discrete delta functions are used to approximate the beam diameter. Such an approach can be found in~\cite{Tschisgale2020, Peskin2006, wu2019, wang2019, kerfriden2020}.
The image on the right shows a representative mesh as commonly applied for the fluid-beam interaction (FBI) approach first proposed for a multi-physics problem in our previous work~\cite{hagmeyer2022}, which serves as the starting point for the present work. The surrogate model is based on the assumption of small beam radii compared to the background elements and couples the two fields only along the beams' centerline instead of its surface. This approximation leads to the possibility of relatively coarse and, thus, computationally efficient background meshes. This in turn makes this method one of the first to be able to approximate the actual physical interaction of the slender bodies while maintaining the complexity reduction with regard to the degrees of freedom gained through the employment of one-dimensional beam theory.

Nevertheless, the solution of mixed-dimensional fluid-structure interaction problems leads to additional algorithmic difficulties. Because of their considerable slenderness and, thus, their susceptibility to external forces, immersing beams within fluid flow generally leads to numerical challenges in the context of the well-known added mass effect~\cite{CAUSIN2005, gerbeau2003, Bukac2016}. Strong coupling of the fluid and structure field is thus crucial for the stability of the simulation. In general, strong coupling can be achieved by solving the full monolithic multi-physics problem, or by means of a staggered partitioned approach with a sufficiently accurate coupling tolerance. So far, the monolithic solution of the arising systems of linearized equations in mixed-dimensional modeling poses many open questions w.r.t. computational efficiency and scalability on parallel computing clusters. In particular, the beam field usually results in non-diagonally dominant (sub-)matrices leading to the need for advanced preconditioning. In addition, the mixed-dimensional nature of the problem calls for specially-tailored preconditioning for the overall monolithic problem as analyzed in~\cite{kuchta2019, budisa2020, hodneland2021}. These challenges and open questions in the context of efficient linear solvers can be bypassed by choosing a strongly-coupled staggered partitioned solution approach. In that case,  well established solution techniques can be applied for the two single fields, that are solved independently from each other. This approach comes with the need for efficient acceleration techniques to guarantee convergence of the staggered scheme because of the added mass effect~\cite{uekermann2021, kuttler2008fixed, gerbeau2003, kiendl2020, Henshaw2015}.

A further algorithmic difficulty stems from the fact that general applications of the envisioned FBI problems exhibit large displacements, making an immersed boundary type method (IBM) the methodology of choice. Such methods have also been used for mixed-dimensional methods in~\cite{wu2019, wang2019, Tschisgale2020, kuchta2021, kerfriden2020}.
In contrast to classical fluid-structure interaction (FSI) approaches using an Arbitrary Lagrangean Eulerian (ALE) description~\cite{mayr2014, Kloppel2012, richter2010, wick2011},
IBM type approaches do not track the interaction interface,
but capture its position as the beams move through the fluid's background mesh.
Once the intersection of a beam with the background mesh has been identified, coupling quantities have to be exchanged between the beam mesh and the fluid's background mesh via mesh tying approaches. Mesh tying is part of many numerical applications ranging from mesh partitioning approaches in the context of high performance computing (HPC)~\cite{wohlmuth2000} to the simulation of multi-physics problems such as FSI~\cite{Kloppel2012, mayr2014, Mayr2020a}. In particular for such surface-coupled problems, the mortar finite element method represents a mathematically sound and stable approach to transfer coupling quantities from one mesh to the other. In that case, the mortar method ensures a straightforward choice of shape functions for the Lagrange multiplier field resulting in a stable mixed finite element formulation for the discretization of the underlying problem. Successful application examples include the case of contact mechanics~\cite{popp2009, popp2010} as well as ALE based fluid-structure interaction~\cite{Kloppel2012, mayr2014, Baaijens2001}.
In the case of surface-to-volume coupling, as applied in immersed boundary type methods, this simple choice of shape functions for inf-sup stable Lagrange multipliers is not guaranteed anymore. Nevertheless, such problems can still benefit in terms of smoothness and robustness of the solution from using a mortar-type segment-to-segment coupling approach. Examples of this are illustrated for fictitious domain/Lagrange multiplier methods in the context of fluid-structure interaction in~\cite{Baaijens2001, santo2008, gerbeau2022} and its advantages have recently been demonstrated for a mixed-dimensional mesh tying problem of 1D beam centerlines embedded into a 3D solid volume~\cite{steinbrecher2020}.
This motivates the use of such a mortar-type discretization approach for the simulation of slender structures interacting with incompressible fluid flow.

In this contribution, the slender structures are modeled using geometrically exact beam theory and embedded into the incompressible Navier-Stokes equations in three dimensions, thus leading to a mixed-dimensional multi-physics problem. The coupling is achieved on the one-dimensional beam centerline using a mortar-type regularized Lagrange multiplier method, which is referred to as a mortar-penalty approach in the following. We apply a Quasi-Newton based Dirichlet-Neumann-type algorithm in combination with a computationally parallel interface capturing technique to the arising mixed-dimensional system. To the best of the authors' knowledge, this is the first time that the algorithmic aspects of a computational framework for the solution of a truly 1D-3D coupled fluid-beam interaction problem have been reported.

The remainder of this paper will be structured as follows: In Section \ref{sec:equations}, we will recount the single field equations as well as the discretization methods used to treat the partial differential equations (PDEs). Section \ref{sec:numerics} will present the continuous Lagrange multiplier based FBI problem, its spatial discretization, the employed penalty-based regularization technique for the Lagrange multiplier degrees of freedom, and the arising discretized nonlinear system of equations. In Section \ref{sec:algo}, we will introduce the algorithmic building blocks for the employed solution algorithm, namely the Dirichlet-Neumann algorithm for the penalty-coupled, immersed method and the acceleration techniques as well as their modifications to suit the proposed approach. In Section \ref{sec:examples}, we will proceed to illustrate the performance of the introduced algorithmic building blocks, and demonstrate the well-suitedness of the proposed mixed-dimensional mortar-type approach by studying its behavior under uniform mesh refinement with respect to a 3D reference solution. We further investigate the influence of the penalty parameter and present an application-based example problem of a submerged vegetation patch showcasing the method's immense potential for systems with a myriad of embedded fibers. Section \ref{sec:conclusions} will provide concluding remarks.

\section{Single Field Equations}
\label{sec:equations}
In this section, we briefly recapitulate the single field equations for the incompressible Navier-Stokes equations as well as the applied beam formulation and recount the employed discrete systems of equations arising from their successive discretization, including the necessary stabilization terms.
We will use the indices $\indBeam$ and $\indFluid$, or  $\mathcal{B}$ and $\mathcal{F}$, to highlight beam- and fluid-related quantities, respectively.
The discretization in space will be solely based on the finite element method (FEM). For the discretization in time, a finite differencing scheme is applied.
In the following, we will assume that the finite differencing scheme can be represented as a linear combination of the state during two successive time steps $n$ and $n+1$, respectively. 
Field-specific time integration schemes will advance from the initial time~$t=0$ until the end time~$T$ of the simulation.
For the sake of a compact notation, we will thus only write the discrete system of equations for the state of one time step $\left[\cdot\right]_{n+1}$. However, this does not constitute a restriction on the time discretization methods itself and will be further addressed in Remark \ref{remark:beam_timint}.

\subsection{Incompressible Navier-Stokes Equations}
\label{subsec:nse}

To model the fluid field, the instationary, incompressible Navier-Stokes equations for Newtonian fluids on fixed meshes are used.

In the following, $H_0^1\left(\Omega_{\indFluid}\right)$ represents the standard Sobolev space on the fluid domain $\Omega_{\indFluid}$ with zero trace on the boundary $\partial \Omega_{\indFluid}$, and we will frequently use the notation $\left(\cdot , \cdot \right)_{X}:=\left(\cdot , \cdot \right)_{L^2\left(X\right)}$ to denote the $L^2$ inner product on the domain $X$. The fluid boundary $\partial\Omega_{\indFluid} = \Gamma^{\indFluid}_N \cup \Gamma^{\indFluid}_D$ can be partitioned into a Neumann boundary $\Gamma^{\indFluid}_N$, on which a traction $\mathbf{h}^{\indFluid}$ is imposed, and a Dirichlet boundary $\Gamma^{\indFluid}_D$, for which a function $\mathbf{v}^{\indFluid}_D$ with prescribed velocity values on $\Gamma^{\indFluid}_D$ is introduced. Furthermore, let $\mathcal{V}_p = \left\{p\in L^2\left(\Omega_{\indFluid}\right) | \left\|p\right\|_{L^2\left(\Omega^{\indFluid}\right)} = 0 \right\}$ be the space of normalized pressure solutions.

For the fluid velocity $\mathbf{v}^{\indFluid}$ and pressure $p^{\indFluid}$, respectively, and using the test functions $\delta\mathbf{v}^{\indFluid}$ and $\delta p^{\indFluid}$, we introduce the semi-linear form
\begin{subequations}
\begin{equation}
	\begin{split}
		\mathbf{a}^{\indFluid}\left(\mathbf{v}^{\indFluid}, p^{\indFluid}; \delta\mathbf{v}^{\indFluid}, \delta p^{\indFluid}\right) &:= \rho^{\indFluid} \left(\dfrac{\partial \mathbf{v}^{\indFluid}}{\partial t}, \delta \mathbf{v}^{\indFluid}\right)_{\Omega_{\indFluid}} + 2\gamma^{\indFluid} \left(\boldsymbol{\mathcal{E}}\left(\mathbf{v}^{\indFluid}\right), \boldsymbol{\nabla} \mathbf{v}^{\indFluid}\right)_{\Omega_{\indFluid}} - \left(p^{\indFluid}, \boldsymbol{\nabla} \cdot \delta\mathbf{v}^{\indFluid}\right)_{\Omega_{\indFluid}} \\
		&+ \rho^{\indFluid}\left(\left(\mathbf{v}^{\indFluid}\cdot \boldsymbol{\nabla}\right)\mathbf{v}^{\indFluid}, \delta \mathbf{v}^{\indFluid}\right)_{\Omega_{\indFluid}} + \left(\boldsymbol{\nabla}\cdot \mathbf{v}^{\indFluid}, \delta p^{\indFluid}\right)_{\Omega_{\indFluid}},
	\end{split}
\end{equation}
and the linear form
\begin{equation}
	\label{eq:force_nse}
	\mathbf{b}^{\indFluid}\left(\delta \mathbf{v}^{\indFluid}\right):=\rho^{\indFluid}\left(\mathbf{f}^{\indFluid}, \delta \mathbf{v}^{\indFluid}\right)_{\Omega_{\indFluid}} + \left(\mathbf{h}^{\indFluid}, \delta \mathbf{v}^{\indFluid}\right)_{\Gamma_N^{\indFluid}},
\end{equation}
\end{subequations}
with the strain rate tensor $\boldsymbol{\mathcal{E}}\left(\mathbf{v}^{\indFluid}\right)$, a body force $\mathbf{f}^{\indFluid}$, the dynamic viscosity $\gamma^{\indFluid}$, and the fluid density $\rho^{\indFluid}$, respectively.

In the case of a divergence-free initial velocity field $\mathbf{v}_0^{\indFluid}$, the behavior of the fluid on the domain $\Omega_{\indFluid}$ is fully described by:

\begin{problem}
	\label{prob:nse}
	Find $\left(\mathbf{v}^{\indFluid}, p^{\indFluid}\right)\in L^2\left(I, H_0^1\left(\Omega_{\indFluid}\right)^3+\mathbf{v}_D\right)
	\times L^2\left(I, \mathcal{V}_p \right)$, with $\mathbf{v}^{\indFluid}=\mathbf{v}_0$ a. e. for $t=0$, such that
	\begin{equation*}
		\int\limits_0^{\trans} \mathbf{a}^{\indFluid}\left(\mathbf{v}^{\indFluid}, p^{\indFluid}; \delta\mathbf{v}^{\indFluid}, \delta p^{\indFluid}\right)
		- \mathbf{b}^{\indFluid}\left(\delta \mathbf{v}^{\indFluid}\right) \text{ d}t=0 \quad \forall \quad \left(\delta \mathbf{v}^{\indFluid}, \delta p^{\indFluid}\right) \in H_0^1\left(\Omega_{\indFluid}\right)^3\times \mathcal{V}_p.
	\end{equation*}
\end{problem}

The discretization of the Navier-Stokes equations leads to a mixed finite element formulation in the fluid velocity and pressure, which generally has to adhere to the inf-sup condition to produce unique and stable solutions. Instead of the employment of inf-sup stable finite elements, we choose to circumvent the condition by adding stabilization terms to the system. This makes it possible to use equal-order P1/P1 shape functions $N_i$ to define the finite element approximations $\fvelfunch$ and $\fpresfunch$ to the fluid solution pair:
\begin{equation}
	\label{eq:fluid_trialfunctions}
	\fvelfunch:=\sum_{i=1}^{n_v^{\indFluid}}N_i\fvelnod, \quad \fpresfunch:= \sum_{i=1}^{n_p^{\indFluid}}N_i\fpresnod,
\end{equation}
Here, $\fvelnod$ and $\fpresnod$ represent the nodal values for the fluid velocity and pressure, respectively, and $n_v^{\indFluid}$ and $n_p^{\indFluid}$ denote the number of nodes for the fluid velocity and fluid pressure, respectively. Afterwards, a PSPG stabilization~\cite{pspg1991} is applied to the Navier-Stokes equations to circumvent the inf-sup condition. To further improve the numerical properties of the resulting system of equations, a div-grad stabilization term is added, and the SUPG method aims at preventing instabilities in convection-dominated flows~\cite{hansbo1990, franca1992, schott2015}.
We can specify the additional stabilization terms for the inf-grad, and SUPG/PSPG term, respectively, by defining the discrete residuals, $\mathbf{r}_h^M$ for the momentum, and $\mathbf{r}_h^C$ for the continuity equation, as
\begin{equation*}
	\mathbf{r}_h^M=\rho^{\indFluid} \dfrac{\partial \fvelfunch}{\partial t} + 2\gamma^{\indFluid} \boldsymbol{\mathcal{E}}\left(\fvelfunch\right) - \boldsymbol{\nabla}\fpresfunch + \rho^{\indFluid}\left(\fvelfunch\cdot \boldsymbol{\nabla}\right)\fvelfunch,
\end{equation*} 
and 
\begin{equation*}
	\mathbf{r}_h^C=\boldsymbol{\nabla}\cdot\fvelfunch.
\end{equation*}
With these definitions, the overall stabilized (time-continuous) nonlinear fluid system reads
\begin{equation}
	\label{eq:nse_semidiscform}
	\begin{split}
		\mathbf{a}_h^{\indFluid}\left(\fvelfunch, \fpresfunch, \delta\fvelfunch, \delta \fpresfunch\right)&:= \rho^{\indFluid} \left(\dfrac{\partial \fvelfunch}{\partial t}, \delta \fvelfunch\right)_{\Omega^h_{\indFluid}} + 2\gamma^{\indFluid} \left(\boldsymbol{\mathcal{E}}\left(\fvelfunch\right), \boldsymbol{\nabla} \fvelfunch\right)_{\Omega^h_{\indFluid}} \\
		&- \left(\fpresfunch, \boldsymbol{\nabla} \cdot \delta\fvelfunch\right)_{\Omega^h_{\indFluid}}
		+ \rho^{\indFluid}\left(\left(\fvelfunch\cdot \boldsymbol{\nabla}\right)\fvelfunch, \delta \fvelfunch\right)_{\Omega^h_{\indFluid}} + \left(\boldsymbol{\nabla}\cdot \fvelfunch, \delta \fpresfunch\right)_{\Omega^h_{\indFluid}} \\
		&+ \rho^{\indFluid}\left(\tau_M \mathbf{r}_h^M, \fvelfunch\cdot\delta\fvelfunch\right)_{\Omega^h_{\indFluid}} + \left(\tau_C\mathbf{r}_{h}^C, \boldsymbol{\nabla}\cdot\delta\fvelfunch\right)_{\Omega^h_{\indFluid}} + \left(\tau_M\mathbf{r}_h^M,\mathbf{\nabla}\delta \fpresfunch\right)_{\Omega^h_{\indFluid}}.
	\end{split}
\end{equation}
Here, $\tau_C$ and $\tau_M$ represent suitable stabilization parameters.
Evaluation of the integrals over the trial and test functions in~\eqref{eq:nse_semidiscform} leads to the nonlinear matrix-vector equations
\begin{equation*}
	\boldsymbol{A}_\mathcal{FF}\left(\mathbf{\hat{v}}_h^{\indFluid}\right)\left(\begin{array}{c}
		\fvelvec \\
		\fpresvec
	\end{array}\right) = \mathbf{f}^\mathcal{F},
\end{equation*}
where the vectors $\fvelvec$ and $\fpresvec$ are made up of the nodal components $\fvelnod$ and $\fpresnod$, respectively, as defined in~\eqref{eq:fluid_trialfunctions}, $\boldsymbol{A}_\mathcal{FF}\left(\fvelvec\right)$ is obtained by evaluation of~\eqref{eq:nse_semidiscform}, and $\mathbf{f}^\mathcal{F}$ contains all nodal contributions of the applied forces introduced in~\eqref{eq:force_nse}.

\subsection{Geometrically Exact Torsionfree Beam Model}
\label{subsubsec:tf}

\begin{figure}
	\begin{center}
		\resizebox{0.4\textwidth}{!}{
			\begin{tikzpicture}
				\draw[thick] (-7, 0)  -- (-2, 0);
				\draw[thick] (-7, 0.2)  -- (-7, -0.2) node[below] {\huge{$0$}};
				\draw[thick] (-2, 0.2)  -- (-2, -0.2) node[below] {\huge{$l$}};
				\draw[thick] (-4.5, 0.2) to[out=90,in=180] (1,3.5);
				\draw[thick] (1,3.5) to (0.7, 3.8);
				\draw[thick] (1,3.5) to (0.7, 3.2);
				\node at (-4, 3.5) {\huge{$\mathbf{r}\left(t, \cdot\right)$}};
				\pgfmathsetseed{2030}
				\draw (0, 3/2) pic[rotate=-40] {bent cylinder};
				\draw (1.5,3.5) -- (3.5,3.5) node[anchor=west] {\huge{$\Omega_{\indBeam}$}};
			\end{tikzpicture}
		}
		\caption{Depiction of a beam and its centerline representation by a curve}
		\label{fig:centerline}
	\end{center}
\end{figure}
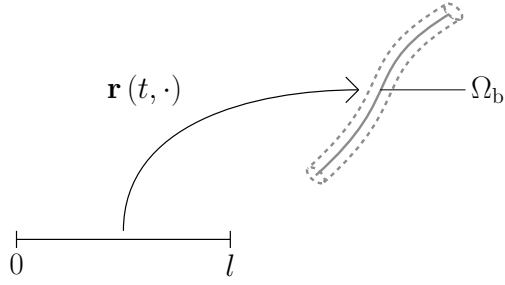
In the case of initially straight beams with isotropic cross-section under the assumption of vanishing shear strains, and assuming that torsional components of distributed and discrete external moments acting on the beam are negligible, the geometrically exact beam theory dating back to Simo~\cite{simo1985} and Reissner~\cite{Reissner1972} can be simplified to a torsionfree (TF) model as proposed in~\cite{Meier2015, meier_contact}. The TF formulation is tailored to high slenderness ratios of the beams and has been found to aid in the avoidance of locking effects in cases where the model requirements are fulfilled. This makes it a desirable choice for the use within the FBI method, which also specifically targets highly slender immersed beams.

\begin{remark}
\label{rem:BeamModels}
We stress that the FBI approach is in no way limited to the TF beam formulation but can straightforwardly be applied to more general Simo-Reissner beam formulations also accounting for torsion as illustrated in the authors' own work~\cite{hagmeyer2022} for the case of one-way coupled problems. Nonetheless, the assumptions summarized above for the TF formulation are approximately fulfilled for the numerical experiments presented in Section~\ref{sec:examples}. For the sake of brevity, we therefore restrict all derivations to the case of the TF formulation.
\end{remark}

In the TF formulation, the state of a beam at each time $t\in I:=\left[0,T\right]$ can be characterized by the position of its centerline $\boldsymbol{r}\left(t, \cdot\right)\in \mathcal{V}_r := \left\{v: \left[0,l\right] \rightarrow \mathbb{R}^3:\left\|\dfrac{\partial^\alpha v}{\partial s}\right\|_\mathcal{B}<\infty, \alpha \leq 2 \right\}$. Figure~\ref{fig:centerline} visualizes the image of the centerline curve $\Omega_{\indBeam}:=\Omega_{\indBeam}\left(t\right):=\mathbf{r}\left(t, \left[0,l\right]\right)\subset\mathbb{R}^3$ denoting the beam's current configuration for any given time $t$.

As is common practice in beam theory, in the following, we will use the abbreviations $\dot{x}:=\dfrac{\partial x}{\partial t}$  and $\ddot{x}:=\dfrac{\partial^2 x}{\partial^2 t}$ for the first and second time derivatives of $x$ as well as ${x}^\prime:=\dfrac{\partial x}{\partial s}$ and ${x}^{\prime\prime}:=\dfrac{\partial^2 x}{\partial^2 s}$ for the first and second derivative with respect to the curve parameter $s$.

Given the density $\rho^{\indBeam}$, the Young's modulus $E$, the length~$l$, the cross-section area $A$ and the torsional moment of inertia $I$, the semi-linear form for the TF model reads
\begin{subequations}
\begin{equation}
	\begin{split}
		\mathbf{a}^{\indBeam}\left(\boldsymbol{r}; \delta\boldsymbol{r}\right) & := \left(\rho^{\indBeam} A\mathbf{\mathbf{\ddot{r}}}, \delta\mathbf{r}\right)_\mathcal{B} + \left(EA \dfrac{\left(\left\|\mathbf{r}^\prime\right\|-1\right)\mathbf{r}^\prime}{\left\|\mathbf{r}^\prime\right\|}, \delta\mathbf{r}^\prime\right)_\mathcal{B}\\ &+ \left(\dfrac{EI}{\left\|\mathbf{r}^\prime\right\|^4}\mathbf{S}\left(\mathbf{r}^\prime\right)\mathbf{r}^{\prime\prime}, \mathbf{S}\left(\delta \mathbf{r}^\prime\right)\mathbf{r}^{\prime\prime}+\mathbf{S}\left(\mathbf{r}^\prime\right)\delta \mathbf{r}^{\prime\prime}\right)_\mathcal{B} \\
		&- \left(\dfrac{2EI}{\left\|\mathbf{r}^\prime\right\|^6}\mathbf{S}\left(\mathbf{r}^\prime\right)\mathbf{r}^{\prime\prime}, \delta \left(\mathbf{r}^{\prime T}  \mathbf{r}\right) \mathbf{S}\left(\mathbf{r}^\prime\right)\mathbf{r}^{\prime\prime}\right)_\mathcal{B} - 
		\left(\mathbf{m},\dfrac{\mathbf{S}\left(\mathbf{r}^\prime\right)}{\left\|\mathbf{r}^\prime\right\|^2} \delta\mathbf{r}^\prime\right)_\mathcal{B} \\
		&- \left[\dfrac{\tilde{\mathbf{m}}\mathbf{S}\left(\mathbf{r}^\prime\right)\delta \mathbf{r}^\prime}{\left\|\mathbf{r}^\prime\right\|^2}\right|_0^l,
	\end{split}
\end{equation}
and the linear form can be given as
\begin{equation}
	\begin{split}
		\mathbf{b}^{\indBeam}\left(\delta\boldsymbol{r}\right):=\left(\mathbf{f}, \delta \mathbf{r}\right)_\mathcal{B}+\left[\mathbf{\tilde{f}}^{\trans}, \delta \mathbf{r}\right|_0^l.
	\end{split}
\end{equation}
\end{subequations}
Here, we denote the cross product of two vectors $\mathbf a$ and $\mathbf b$ by the cross product operator $\boldsymbol{S}\left(\boldsymbol{a}\right)\boldsymbol{b}:=\mathbf{a} \times \mathbf{b}$ and use the external forces $\mathbf{f}$, the external moments $\mathbf{m}$, the point forces $\mathbf{\tilde{f}}$, and the point moments $\mathbf{\tilde{m}}$.

Now, we can complete the general beam problem by introducing the Dirichlet boundary conditions $\boldsymbol{r}_D$ and an initial state $\boldsymbol{r}_0$. Then, the state of the beam can be characterized as solution $\boldsymbol{r}$ of the following problem:

\begin{problem}
	\label{prob:beam}
	Find $\boldsymbol{r}\in L^2\left(I, \mathcal{V}_r+\boldsymbol{r}_D\right)$ such that
	\begin{equation}
		\label{eq:beamproblem}
		\int\limits_0^{\trans} \mathbf{a}^{\indBeam}\left(\boldsymbol{r}; \delta\boldsymbol{r}\right) - \mathbf{b}^{\indBeam}\left( \delta\boldsymbol{r}\right) \text{ d}t=0 \qquad \forall \quad \delta\boldsymbol{r}\in \mathcal{V}_r,
	\end{equation}
	with $\boldsymbol{r}=\boldsymbol{r}_0$ a. e. for $t=0$, and $\left\|\mathbf{r}^\prime\left(0, \cdot\right)\right\|_{\mathbb{R}^3}=1$ a. e. on $\left[0,l\right]$.
\end{problem}

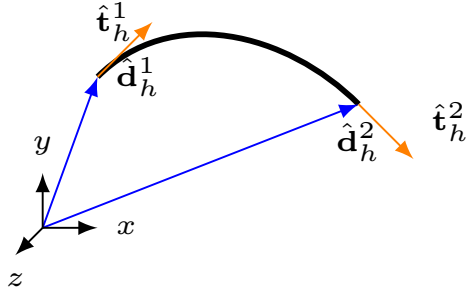
\begin{figure}
	\begin{center}
	\resizebox{0.4\textwidth}{!}{
\begin{tikzpicture}

\draw[very thick] (-1.6, 0) to [out=45, in=135] (0.3, -0.2);

\draw[-latex, orange] (-1.6, 0) -- (-1.2, 0.4) node[anchor=east]{\color{black}\tiny{$\hat{\mathbf{t}}_{h}^1$}};
\draw[-latex, orange] (0.3, -0.2) -- (0.7, -0.6) node[anchor=south west]{\color{black}\tiny{$\hat{\mathbf{t}}_{h}^2$}}; 

\draw[-latex, blue] (-2, -1.1) -- (-1.6, 0) node[anchor=west]{\color{black}\tiny{$\hat{\mathbf{d}}_{h}^1$}};
\draw[-latex, blue] (-2, -1.1) -- (0.3, -0.2) node[anchor=north]{\color{black}\tiny{$\hat{\mathbf{d}}_{h}^2$}};

\draw[-latex] (-2, -1.1) -- (-1.6, -1.1) node[anchor=west]{\tiny{$x$}};
\draw[-latex] (-2, -1.1) -- (-2, -0.7) node[anchor=south]{\tiny{$y$}};
\draw[-latex] (-2, -1.1) -- (-2.2, -1.3) node[anchor=north]{\tiny{$z$}};

\end{tikzpicture}}
	\caption{Illustration of the beam finite element interpolation using 3rd-order $C^1$-continuous Hermite shape functions}
	\label{fig:beamelement}
\end{center}
\end{figure}

In this work, the beam centerline and velocity discretization is exclusively based on 3rd-order $C^1$-continuous Hermite finite elements
as sketched in Figure~\ref{fig:beamelement}.
For detailed information on the construction of these shape functions, the interested reader is referred to~\cite{Meier2015,meier_contact}.
For now, let it be said that $C^1$-continuity is required to capture the beam's curvature and, thus, for the torsionfree beam formulation presented above to be well-defined.
Based on the nodal positions~$\rnodtrans$, the nodal tangents~$\rnodtang$, and the corresponding scalar-valued positional and tangential shape functions~$H^d_j$ and~$H^t_j$, respectively,
the semi-discrete centerline position and beam velocity field can be expressed as
\begin{equation}
	\label{eq:beam_fe}
	\rfunch:=\sum_{j=1}^{n^{\indBeam}}H^d_j\rnodtrans + \dfrac{l_{ele}}{2} H^t_j\rnodtang, \quad  \bvelfunch:=\dot{\mathbf{r}}_h.
\end{equation}
Here, $n^{\indBeam}$ represents the number of beam nodes.
Explicit formulas for the shape functions~$H^d_j$ and~$H^t_j$ are given as part of Appendix~\ref{app:ExampleMortarMatrices}.

Introduction of~\eqref{eq:beam_fe} into~\eqref{eq:beamproblem} and subsequent spatial integration leads to the nonlinear system 
\begin{equation*}
	\mathbf{A}_\mathcal{BB}\left(\rvec\right)\rvec = \brhsvec.
\end{equation*}
Here, $\rvec$ contains all beam-related degrees of freedom (DoFs) including positional and tangential unknowns, and accordingly so do the matrix $\mathbf{A}_\mathcal{BB}$ and the vector $\brhsvec$.

\begin{remark}
	\label{remark:beam_timint}
	In our case, the Generalized-$\alpha$ time integration scheme~\cite{genalpha} to approximate the beams' velocity and acceleration does not follow the assumption of representing a linear combination of the solutions of two successive time steps. We use the equations
	\begin{subequations}
	\begin{eqnarray}
		\label{eq:timeintegration_beam}
		\dot{\boldsymbol{r}}_{n+1} = \dfrac{\gamma}{\beta \Delta t}\left(\boldsymbol{r}_{n+1} - \boldsymbol{r}_{n}\right) - \dfrac{\gamma-\beta}{\beta} \dot{\boldsymbol{r}}_{n} - \dfrac{\gamma - 2\beta}{2\beta}\Delta t\ddot{\boldsymbol{r}}_{n}, \\
		\ddot{\boldsymbol{r}}_{n+1} = \dfrac{1}{\beta \delta t^2}\left(\boldsymbol{r}_{n+1}-\boldsymbol{r}_{n}\right) - \dfrac{1}{\beta \Delta t}\dot{\boldsymbol{r}}_{n} - \dfrac{1-2\beta}{2\beta}\ddot{\boldsymbol{r}}_{n}, \\
		\mathbf{G}_{\text{int}, n+1-\alpha_f} + \mathbf{G}_{\text{kin}}\left(\ddot{\boldsymbol{r}}_{n+1-\alpha_m}, \delta\mathbf{r}\right)- \mathbf{G}_{\text{ext}, n+1-\alpha_f} = 0,
	\end{eqnarray}
	\end{subequations}
	where we compute the internal energy $\mathbf{G}_{\text{int}, n+1-\alpha_f}$ and the kinetic energy $\mathbf{G}_{\text{kin}}\left(\ddot{\boldsymbol{r}}_{n+1-\alpha_m}\right)$ using the trapezoidal rule instead of using a linear combination of $\mathbf{G}_{\text{int}, n+1}$ and $\mathbf{G}_{\text{int}, n}$, or $\mathbf{G}_{\text{kin}}\left(\ddot{\boldsymbol{r}}_{n+1}\right)$ and $\mathbf{G}_{\text{kin}}\left(\ddot{\boldsymbol{r}}_{n}\right)$, respectively. For the sake of a compact notation, we chose to neglect this fact within this work. Nevertheless, all methods and discretization schemes presented can still be applied in the same manner. Note also, that the aforementioned method to compute $\mathbf{G}_{\text{int}, n+1-\alpha_f}$ and $\mathbf{G}_{\text{kin}}\left(\ddot{\boldsymbol{r}}_{n+1-\alpha_m}\right)$, namely the trapezoidal and midpoint rule, are not equal here as we are using nonlinear finite element methods. However, both methods coincide for linear finite element technology.
\end{remark}

\section{Fluid Beam Interaction Problem}
\label{sec:numerics}
We will first state the full continuous FBI interaction problem by defining a Lagrange multiplier field $\lambda$ along the beam centerline to enforce the coupling constraints.
Afterwards, we introduce a mortar-type finite element discretization for this Lagrange multiplier field as well as a successive regularization thereof. Finally, we derive all necessary coupling operators, and present the overall coupled nonlinear fluid-beam interaction system.

\subsection{Continuous Fluid-Beam Interaction Constraints}
\label{sec:GovEqFBICoupling}

As already discussed in Section \ref{sec:intro}, the aim of the applied FBI method is to yield a computationally efficient approach that does not require the geometrical resolution of the embedded fibers by the background mesh as in~\cite{Tschisgale2020, wang2019, Henshaw2015, tschisgale2021}. Instead, we assume the immersed fibers to be highly slender, also compared to the background mesh, and follow the idea in~\cite{Baaijens2001, steinbrecher2020, wu2019, kremheller2019, Peskin2002} of coupling the fluid to the fiber directly on the one-dimensional beam centerline leading to an efficient surrogate model for problems involving immersed slender structures. A comprehensive discussion can be found in the authors' previous contribution~\cite{hagmeyer2022}. Using a Lagrange multiplier approach to enforce the kinematic coupling condition, i.e. the continuity of beam and fluid velocities $\mathbf{v}^{\indBeam}$ and $\mathbf{v}^{\indFluid}$, respectively, we obtain the overall FBI problem:

\begin{problem}
	Find $\left(\mathbf{v}^{\indFluid}, p^{\indFluid}, \boldsymbol{r}, \boldsymbol{\lambda}\right)\in L^2\left(I, H^1\left(\Omega_{\indFluid}\right)^3\right)\times L^2\left(I, \mathcal{V}_p \right)\times L^2\left(I, \mathcal{V}_r\right)\times L^2\left(I, \mathcal{V}_r\right)$, with $\mathbf{v}^{\indFluid}=\mathbf{v}_0$, $\boldsymbol{r}=\boldsymbol{r}_0$ a. e. for $t=0$, and $\left\|\mathbf{r}^\prime\left(0, \cdot\right)\right\|_{\mathbb{R}^3}=1$, such that
	\begin{equation}
		\label{eq:fbi}
		\begin{split}
			&\int\limits_0^T \mathbf{a}^{\indFluid}\left(\mathbf{v}^{\indFluid}, p^{\indFluid}; \delta\mathbf{v}^{\indFluid}, \delta p^{\indFluid}\right)
			- \mathbf{b}^{\indFluid}\left(\delta \mathbf{v}^{\indFluid}\right) + \int\limits_0^l \boldsymbol{\lambda}\boldsymbol{\Pi}\delta\mathbf{v}^{\indFluid}\circ \mathbf{r} \text{d}s \text{ d}t = 0, \\
			&\int\limits_0^T \mathbf{a}^{\indBeam}\left(\boldsymbol{r}; \delta\boldsymbol{r}\right) - \mathbf{b}^{\indBeam}\left(\delta\boldsymbol{r}\right) - \int\limits_0^l \boldsymbol{\lambda}\delta\mathbf{r} \text{d}s \text{ d}t =0, \\
			&\int\limits_0^T\int\limits_0^l\left(\mathbf{v}^{\indFluid}-\mathbf{v}^{\indBeam}\right)\delta\boldsymbol{\lambda}\text{d}s \text{ d}t = 0
		\end{split}
	\end{equation}
	for all $\left(\delta\mathbf{v}^{\indFluid}, \delta p^{\indFluid}, \delta \boldsymbol{r}, \delta \boldsymbol{\lambda}\right)\in H_0^1\left(\Omega_{\indFluid}\right)^3\times \mathcal{V}_p\times \mathcal{V}_r\times \mathcal{V}_r$.
	\label{prob:fbi}
\end{problem}
Here, $\boldsymbol{\lambda}$ is the Lagrange multiplier, and $\boldsymbol{\Pi}$ represents the 1D-3D analog to the projection postulated by the well known trace theorem,
as the fluid variable $\mathbf{v}^{\indFluid}$ is a priori not well defined on $[0,l]$.
For further discussions on the nature of the projection operator $\boldsymbol{\Pi}$, we refer the interested reader to~\cite{hagmeyer2022, kuchta2021}.
The manuscript at hand, however, focuses on the numerical aspects of the proposed method.
A rigorous mathematical analysis on existence and well-posedness of the continuous problem is out of the scope of this work.

In accordance with the dynamic coupling condition ensuring the continuity of traction, the Lagrange multiplier $\boldsymbol{\lambda}$ can be interpreted as a line load ensuring the equilibrium of forces. Nevertheless, in the same manner as in~\cite{Peskin2002, wong2020, Baaijens2001}, the choice of the coupling domain as the one-dimensional beam centerline effectively introduces a singularity into the continuous system of equations. The handling of this singularity in the numerical system and its effect on the numerical behavior has been demonstrated in the authors' previous work~\cite{hagmeyer2022} and will be shortly discussed in Section \ref{sec:examples} in the context of the method's convergence behavior under uniform mesh refinement.

\begin{remark}
\label{rem:TorsionalCoupling}
Note that rotational effects are assumed to be negligible due to the application to highly slender structures.
While this assumption is valid for a wide range of application scenarios,
there exist cases where rotational effects of a beam on the flow around it and vice versa are not negligible even for high slenderness ratios.
For such scenarios, the presented method needs to be extended in order to include rotational coupling.
We have investigated rotational coupling in the context of the embedding of 1D beams within 3D solid volumes~\cite{steinbrecher2022}
and this necessitates the construction of tailored rotational triads.
In the context of FSI involving slender structures,
coupling of torsional modes has so far only been achieved by re-constructing a two-dimensional fictitious coupling boundary based on the 1D beam centerline~\cite{huang2019}.
For an extension of the present method, the equations and ideas of both these methods will have to be combined.
However, this extension is beyond the scope of this contribution and is subject to future investigations.
To illustrate the wide range of application scenarios for which the above simplification assumptions hold,
Section~\ref{sec:examples} presents a number of academic and application-inspired scenarios.
\end{remark}

\subsection{Regularized Lagrange Multiplier Field and Discretized System of Equations}
\label{sec:RegularizedLagMult}

The introduction of Lagrange multipliers to enforce the continuity of velocities on the coupling interface leads to a mixed finite element formulation.
Assuming Lagrange multiplier shape functions~$\Phi_k$, the Lagrange multiplier field is interpolated from its nodal values via
\begin{equation}
	\label{eq:lm_shape}
	\lambdafunch:=\sum_{k=1}^{n^{\indLagMult}}\Phi_k\lambdanod,
\end{equation}
where $n^{\indLagMult}$ denotes the number of nodes carrying a Lagrange multiplier and the shape functions~$\Phi_k$ need to be chosen carefully to adhere to the inf-sup condition~\cite{Babuska1973a,Brezzi1974a}.
While the choice of such inf-sup stable shape functions for surface-coupled meshes on conforming surfaces can be straightforward,
namely by using the same shape functions as for the primal field and, thus, leading to the mortar finite element method~\cite{wohlmuth2000, popp2009},
constructing the required inf-sup stable finite elements to guarantee well-posedness of problems arising from immersed methods, particularly mixed-dimensional ones, is a complex task
and poses several open research questions.

As an alternative, we circumvent the inf-sup condition by first performing a mortar discretization, but then applying a node-wise penalty regularization,
as for example in contact mechanics~\cite{yang2005penalty, puso2004penalty}.
Thereby, we assume a linear relation between the Lagrange multipliers and the mismatch in velocities,
where the slope of this linear relation is given by the penalty parameter~$\epsilon$:
\begin{equation*}
	\lambdavec :=\epsilon\left(\fvelfunch-\bvelfunch\right)
\end{equation*}
In contrast to classical Gauss-point-to-segment (GPTS) based penalty approaches,
this mortar-penalty regularization approach keeps the Lagrange multiplier shape functions, which now take the role of interpolation functions of the regularized interface variable,
to fulfill the coupling in a weak sense and introduces matrix-matrix products into the nonlinear system.

To represent the arising nonlinear system,
we insert~\eqref{eq:lm_shape} into the coupling constraint~\eqref{eq:fbi} and perform integration along the beam centerline,
yielding the following mortar-type matrices~$\mathbf{D}$ and~$\mathbf{M}$ and the scaling matrix~$\boldsymbol{\kappa}$:
\begin{equation}
\label{eq:mortarmatrices}
\begin{split}
\mathbf{D}\left(\indLagMultDof, \indBeamDof\right)=\int\limits_0^l \Phi_\indLagMultDof H_\indBeamDof\text{ d}s,
\quad
\mathbf{M}\left(\indLagMultDof, \indFluidDof\right)=\int\limits_0^l \Phi_\indLagMultDof N_\indFluidDof\circ \chi\text{ d}s,
\quad
\boldsymbol{\kappa}\left(\indLagMultDof, \indLagMultDof\right)=\int\limits_0^l \Phi_\indLagMultDof\text{ d}s,\\
\text{for } \indLagMultDof = \{1,\hdots, N^{\indLagMult}\},\indBeamDof = \{1,\hdots,N^{\indBeam}\},\indFluidDof = \{1,\hdots,N^{\indFluid}\}
\end{split}
\end{equation}
Using $\ndim$ to refer to the number of spatial dimensions (in our case, we work with $\ndim=3$),
$N^{\indLagMult} := \ndim\cdot n^{\indLagMult}$ denotes the total number of Lagrange multiplier DoFs,
while $N^{\indBeam}:= 2\ndim\cdot n^{\indBeam}$ denotes the total number of beam DoFs.
Assuming an appropriate mapping between beam nodes and their positional and tangential DoFs,
the family of scalar-valued shape functions~$H_\indBeamDof$
summarizes both positional and tangential shape functions~$H^d_i$ and~$H^t_i$ introduced in~\eqref{eq:beam_fe}.
Furthermore, $N^\indFluid:=\ndim\cdot n^\indFluid_v+n^\indFluid_p$ denotes the overall number of unknowns for the fluid field, and $\chi$ represents the projection of the beam centerline onto its current position in three-dimensional space.
Note that this projection introduces a dependency on the current position of the beams to the nonlinear system of equations.
In particular, the sparsity pattern of $\mathbf{M}$ will change as the beams move through the fluid mesh.
The ordering of beam's positional and tangential DoFs can also be chosen in a different form than given above. A concrete example for such an ordering as well as the evaluation of~$\mathbf{D}$ and~$\mathbf{M}$ are given in Appendix~\ref{app:ExampleMortarMatrices}.

Based on the mortar matrices defined in~\eqref{eq:mortarmatrices}, the overall system then takes the following form: 
\begin{equation}
	\label{eq:monolithicsystem}
	\left[\begin{array}{ccc}
		\mathbf{A}_\mathcal{FF} + \epsilon\mathbf{M}^{\trans}\boldsymbol{\kappa}^{-1}\mathbf{M} & 0 & -\epsilon\mathbf{M}^{\trans}\boldsymbol{\kappa}^{-1}\mathbf{D}\\
		-\epsilon\mathbf{D}^{\trans}\boldsymbol{\kappa}^{-1}\mathbf{M} & \mathbf{A}_\mathcal{BB} & \epsilon\mathbf{D}^{\trans}\boldsymbol{\kappa}^{-1}\mathbf{D}
	\end{array}\right]_{n+1}
	\left[\begin{array}{c}
		\left(
		\begin{array}{c}
			\fvelvec \\
			\fpresvec \\
		\end{array}\right)  \\
		\rvec \\
		\bvelvec\left(\rvec\right)
	\end{array}\right]_{n+1}=
	\left[\begin{array}{c}
		\frhsvec \\
		\brhsvec
	\end{array}\right]_{n+1}.
\end{equation}

Here, the normalization with the matrix $\boldsymbol{\kappa}$
is classically necessary to fulfill simple patch tests, and is used to weigh the penalty force exchanged between the two fields~\cite{yang2005penalty}.

\begin{remark}
	For the sake of brevity we use $\mathbf{A}_\mathcal{FF}$ to represent the entire fluid matrix including stabilization and pressure contributions as introduced in Section \ref{subsec:nse}. Since the fluid pressure has no effect on the coupling, the respective terms in the coupling matrix $\mathbf{M}$ thus default to zero. For simplicity, we will further use $\fvelvec:=\left(\begin{array}{c}
		\fvelvec \\
		\fpresvec
	\end{array}\right)$ to represent the fluid velocity as well as fluid pressure degrees of freedom.
\end{remark}

Since the matrix operator~\eqref{eq:monolithicsystem} is rectangular the system in its current form is not solvable. There are two ways to remedy this problem: i) add an additional constraint to the system, or ii) remove one of the unknowns. The obvious target for either of these strategies are the beam velocity unknowns. According to~\eqref{eq:timeintegration_beam}, i) leads to the addition of the line
\begin{equation*}
	\left[
	\begin{array}{ccc}
		0 & -\dfrac{\gamma}{\beta \Delta t} \boldsymbol{1} & \boldsymbol{1} \\
	\end{array}
	\right]_{n+1}
	\left[\begin{array}{c}
		\fvelvec \\
		\rvec \\
		\bvelvec
	\end{array}\right]_{n+1} =
	\left[\begin{array}{c}
		\btrhsvecn
	\end{array}\right]_{n+1},
\end{equation*}
where $\mathbf{f}_{h}^{\mathcal{B}, n}=  -\dfrac{\gamma}{\beta \Delta t}\boldsymbol{r}_{n} - \dfrac{\gamma-\beta}{\beta} \dot{\boldsymbol{r}}_{n} - \dfrac{\gamma - 2\beta}{2\beta}\Delta t\ddot{\boldsymbol{r}}_{n}$. This line represents a weak enforcement of the time integration scheme.

However, we will follow the more commonly used strategy ii), and directly use the time integration scheme to express $\bvelvec\left(\rvec\right)$ in terms of $\rvec$. In that case,~\eqref{eq:monolithicsystem} becomes
\begin{equation}
	\label{eq:monolithicsystem_complete}
		\left[\begin{array}{cc}
			\mathbf{A}_\mathcal{FF} + \epsilon\mathbf{M}^{\trans}\boldsymbol{\kappa}^{-1}\mathbf{M} & - \dfrac{\gamma}{\beta\Delta t}\epsilon\mathbf{M}^{\trans}\boldsymbol{\kappa}^{-1}\mathbf{D}\\
			-\epsilon\mathbf{D}^{\trans}\boldsymbol{\kappa}^{-1}\mathbf{M} & \mathbf{A}_\mathcal{BB} + \dfrac{\gamma}{\beta\Delta t}\epsilon\mathbf{D}^{\trans}\boldsymbol{\kappa}^{-1}\mathbf{D}
		\end{array}\right]_{n+1}
		\left[\begin{array}{c}
			\fvelvec \\
			\rvec
		\end{array}\right]_{n+1}
		=
		\left[\begin{array}{c}
			\frhsvec - \epsilon\mathbf{M}^{\trans}\boldsymbol{\kappa}^{-1}\mathbf{D}\btrhsvecn\\
			\brhsvec + \epsilon\mathbf{D}^{\trans}\boldsymbol{\kappa}^{-1}\mathbf{D}\btrhsvecn
		\end{array}\right]_{n+1}.
\end{equation}
The employed algorithmic strategy to solve this nonlinear system of equations will be discussed in the upcoming Section~\ref{sec:algo}.

\section{Partitioned Solution Algorithm}
\label{sec:algo}

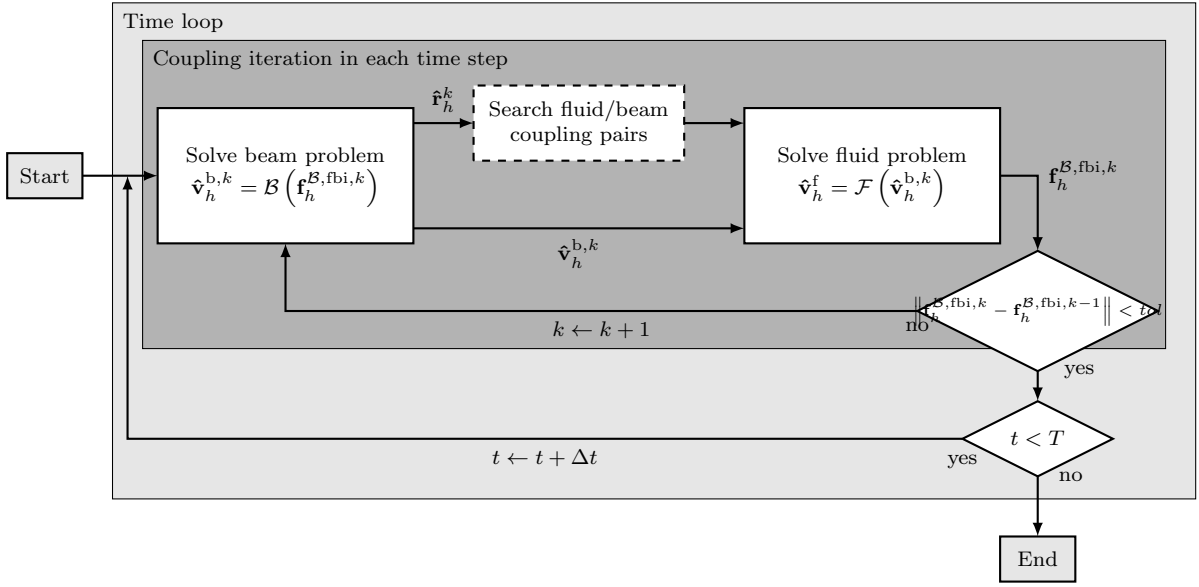
\begin{figure}
\centering

\begin{tikzpicture}[font=\scriptsize]

\tikzstyle{myline} = [-latex,thick]

\draw [fill=gray!20] (-0.2,5.8) node [below right] {Time loop} rectangle (14.2,-0.8);

\draw [fill=gray!60] (0.2,5.3) node [below right] {Coupling iteration in each time step} rectangle (13.8,1.2);

\draw [thick,fill=white] (0.4,2.6) rectangle (3.8,4.4) node [pos=0.5] {\begin{tabular}{c}
Solve beam problem\\
$\bvelveck=\mathcal{B}\left(\bfbiforceveck\right)$
\end{tabular}};

\draw [thick,fill=white] (8.2,2.6) rectangle (11.6,4.4) node [pos=0.5] {\begin{tabular}{c}
Solve fluid problem\\
$\fvelvec=\mathcal{F}\left(\bvelveck\right)$
\end{tabular}};

\draw [myline] (3.8,4.2) -- (4.6,4.2) node [pos=.5,above] {$\rveck$};
\draw [thick,dashed,fill=white] (4.6,3.7) rectangle (7.4,4.7) node [pos=0.5] {\begin{tabular}{c}
Search fluid/beam\\
coupling pairs
\end{tabular}};
\draw [myline] (7.4,4.2) -- (8.2,4.2);

\draw [myline] (3.8,2.8) -- (8.2,2.8) node [pos=.5,below] {$\bvelveck$};

\draw [thick,fill=white] (10.5,1.7) -- (12.1,2.5) -- (13.7,1.7) -- (12.1,0.9) -- cycle;
\node at (12.1,1.7) {\tiny$\left\|\bfbiforceveck-\bfbiforceveckold\right\|<tol$};

\draw [myline] (11.6,3.5) -- (12.1,3.5) node [right] {$\bfbiforceveck$} -- (12.1,2.5);

\draw [myline] (10.5,1.7) node [below] {no} -- (2.1,1.7) node [below,pos=0.5] {$k\leftarrow k+1$} -- (2.1,2.6);

\draw [myline] (12.1,0.9) node [right=6pt] {yes} -- (12.1,0.5);
\draw [thick,fill=white] (11.1,0.0) -- (12.1,0.5) -- (13.1,0.0) -- (12.1,-0.5) -- cycle;
\node at (12.1,0.0) {\scriptsize$t < T$};

\draw [myline] (11.1,0.0) node [below=2pt] {yes} -- (0.0,0.0) node [pos=0.5,below] {$t\leftarrow t+\Delta t$} -- (0.0,3.5);
\draw [myline] (12.1,-0.5) node [right=4pt] {no} -- (12.1,-1.3);

\draw [thick,fill=gray!20] (11.6,-1.3) rectangle (12.6,-1.9) node [pos=0.5] {End};

\draw [thick,fill=gray!20] (-1.6,3.2) rectangle (-0.6,3.8) node [pos=0.5] {Start};
\draw [myline] (-0.6,3.5) -- (0.4,3.5);

\end{tikzpicture}
\caption{Partitioned algorithm to solve the immersed FBI problem by evaluating the problem $\mathcal{F}\left(\bvelveck\right)$
on the fluid partition and the problem $\mathcal{B}\left(\fbiforceveck\right)$ on the beam partition in an iteratively coupled manner.}
\label{fig:algo}
\end{figure}

The partitioned algorithm to solve the interaction problem~\eqref{eq:monolithicsystem} in a staggered manner is visualized in Figure~\ref{fig:algo}.
It has evolved as an extension of mixed-dimensional one-way coupling schemes using a GPTS approach\cite{hagmeyer2022}.
First, the general partitioning and transfer of coupling variables between the fluid and beam partitions will be discussed in Section~\ref{subsec:coupling}.
Acceleration techniques to speed-up the convergence of the partitioned scheme will be introduced in Section~\ref{subsec:accelorator}.
Finally, Section~\ref{sec:NonlinearChallenges} addresses some computational challenges arising from additional nonlinearities,
in particular the dependence of the coupling matrices on the current state as well as the movement of the fibers relative to the background mesh.

\subsection{FBI Coupling}
\label{subsec:coupling}

As visualized in Figure~\ref{fig:algo}, the partitioned algorithm is set-up as a Dirichlet-Neumann-type algorithm where the kinematic coupling condition serves as a Dirichlet condition on the fluid field, and a Neumann boundary condition is applied to the beam partition. Here, the kinematic constraint is enforced weakly on the fluid partition, resembling partitioned algorithms presented in~\cite{dirichletneumannconvergence} for two-body contact problems.
In the following, we will use the subscript $k$ to denote variables computed in the $k$th FBI iteration. For the sake of a compact notation, we will drop this index in equations only involving variables computed within the same FBI iteration. 

To enforce the kinematic constraint for a given beam velocity $\bvelvec$ in a weak sense on the fluid partition, the first line of the coupled system~\eqref{eq:monolithicsystem_complete}, namely
\begin{equation}
	\label{eq:dirichletpartition}
	\left(\mathbf{A}_\mathcal{FF}+\epsilon\mathbf{M}^{\trans}\boldsymbol{\kappa}^{-1}\mathbf{M}\right)\fvelvec = \frhsvec+ \epsilon\mathbf{M}^{\trans}\boldsymbol{\kappa}^{-1}\mathbf{D}\bvelvec,
\end{equation}
is solved for~$\fvelvec$. This nonlinear solution procedure is represented by the solution operator $\mathcal{F}: \bvelvec \rightarrow \fvelvec$ in Figure~\ref{fig:algo}.
In that case, \eqref{eq:dirichletpartition} can be interpreted as the discretization of the Navier-Stokes equations with the addition of a penalization of the constraint error
\begin{equation*}
		\mathbf{e}_{constr} := \mathbf{M}^{\trans}\boldsymbol{\kappa}^{-1}\mathbf{M}\fvelvec - \mathbf{M}^{\trans}\boldsymbol{\kappa}^{-1}\mathbf{D}\bvelvec,
\end{equation*}
effectively enforcing equality of velocities on the coupling interface in a weak sense.
This approach of weak enforcement using a mortar finite element type method has multiple advantages: firstly, in contrast to a strong enforcement, the weak enforcement of the constraint naturally allows for fulfillment of the divergence condition also close to the immersed beams. Secondly, the mortar finite element method provides a natural way for transferring velocity and force data from one mesh to the other. This is of particular interest for the non-matching meshes at hand that unavoidably arise from the dimensionality gap between the coupled 1D and 3D fields.
As a direct result, the interaction force, that is required to keep the equations on the beam partition in equilibrium as it acts on the beam mesh, is stated in the second line of the monolithic system~\eqref{eq:monolithicsystem_complete} as
\begin{equation*}
	\bfbiforcevec = \epsilon\mathbf{D}^{\trans}\boldsymbol{\kappa}^{-1}\mathbf{M}\fvelvec - \epsilon\mathbf{D}^{\trans}\boldsymbol{\kappa}^{-1}\mathbf{D}\bvelvec.
\end{equation*}
For a given fluid velocity $\fvelvec$ and a given beam velocity $\bvelvec$, this force is then applied to the beam partition which leads to the following overall nonlinear system:
\begin{equation}
	\label{eq:beam_fbi_partition}
	\mathbf{A}_\mathcal{BB}\rvec = \brhsvec + \bfbiforcevec
\end{equation}
Solving~\eqref{eq:beam_fbi_partition} for $\rvec$ provides a solution to the primal beam variable, i.e. the beam centerline $\rvec$. However, as the fluid partition expects the beam velocity $\bvelvec:=\bvelvec\left(\rvec\right)$ as input, a suitable time integration scheme $\mathcal{T}$ has to be applied to the beams' current position vector $\rvec$. See Remark~\ref{remark:beam_timint} for an example of such a time integration scheme. The current beam velocity can, therefore, be obtained via 
\begin{equation*}
	\mathcal{B}\left(\bfbiforcevec\right):=\mathcal{T}\left(\mathbf{A}_\mathcal{BB}^{-1}\left( \brhsvec + \bfbiforcevec\right), \rvec\left(t-\Delta t\right),...\right),
\end{equation*}
and can then be passed to the fluid partition for use in the subsequent coupling iteration~$k+1$.

The solution of~\eqref{eq:beam_fbi_partition} generally leads to a new centerline position for the beams.
Especially, the beams' location relative to the fluid background mesh is changed.
As the coupling matrix~$\mathbf{M}$ depends on the projection~$\chi$ of the discrete centerline position onto the fluid mesh,
fluid-beam element pairs have to be identified, a projection to the background elements has to be performed, and the coupling matrices have to be build anew.
As the fluid partition relies on these coupling matrices, it implicitly depends on the vector containing the beam centerline position as well.
We refer to Section~\ref{sec:NonlinearChallenges} for details on tackling the associated computational challenges.

In order to test convergence of the solution scheme, the stopping criterion
\begin{equation*}
	\left\|\bfbiforceveck-\bfbiforceveckold\right\|<tol
\end{equation*}
is applied to the interaction force acting on the beam. Here, $tol$ represents an user-given convergence tolerance.
This completes the Dirichlet-Neumann-type algorithm for immersed beams as sketched in Figure~\ref{fig:algo}.
However, in the case of highly slender fibers, a simple staggered approach as the one described so far may exhibit convergence problems.
The upcoming Section~\ref{subsec:accelorator} addresses this challenge and the choice of a suitable convergence acceleration technique.

\begin{remark}
	\label{remark:segmentation}
	In order to compute the projection $\chi$ introduced in~\eqref{eq:mortarmatrices}, the following constraint has to be fulfilled:
	\begin{equation*}
		\label{eq:segmentation_constraint}
		g\left(\xi_1, \xi_2, \xi_3, \eta\right):= \mathbf{X}\left(\xi_1, \xi_2, \xi_3\right) - \mathbf{R}\left(\eta\right) = 0,
	\end{equation*}
where $\mathbf{X}$ and $\mathbf{R}$ represent the position of the fluid element and beam centerline according to the parameter space and the parameter coordinates $\xi_1, \xi_2, \xi_3$, and  $\eta$, respectively. Solution of~\eqref{eq:segmentation_constraint} for the fluid element's parameters $\xi_1, \xi_2, \xi_3$ finally gives the means to calculate the coupling matrices~\eqref{eq:mortarmatrices}. In general, a local Newton solver is necessary to obtain the solution of~\eqref{eq:segmentation_constraint}. More information on the projection and possible segmentation procedure can be found in~\cite{steinbrecher2020}.
\end{remark}

\begin{remark}
	Here, $\bfbiforcevec$ is held constant within the nonlinear solution procedure of~\eqref{eq:beam_fbi_partition} and is not updated according to $\bvelvec$ during the Newton-Raphson procedure. This nonlinearity, therefore, has to be handled by the FBI algorithm. This, in particular, necessitates the acceleration techniques presented in the upcoming subsection.
\end{remark}

\subsection{Acceleration Technique}
\label{subsec:accelorator}

Due to the partitioning, the applied interaction force remains constant during the solution of the structural problem, which leads to a neglect of the change in geometry of the beam, and therefore, the FBI interface. This nonlinearity is thus not treated within the Newton solver for the beam system, but successively updated through the coupling iterations between the two fields. Especially for very slender and, thus, sensitive structures, this nonlinearity leads to convergence problems of the coupling iterations. Acceleration techniques to facilitate convergence of the partitioned algorithm are well-established for classical 3D FSI problems.
Within this work, we use two different methods to treat the nonlinearity, namely the Aitken relaxation method~\cite{aitken} as a representative of simple fixed-point iterations, and a matrix-free Newton Krylov (MFNK) method~\cite{mfnk2004} as a representative of the somewhat more involved Quasi-Newton methods.

\subsubsection{Aitken Relaxation Method}
\label{subsubsec:aitken}

The idea behind the Aitken method is to relax the exchanged interface variable, classically the interface displacement, in order to avoid oscillations in the convergence behavior as described in~\cite{kuttler2008fixed}. The amount of relaxation is controlled via the relaxation parameter $\omega\in\left(0,1\right)$, which is recomputed based on the current solution within every coupling step. In~\cite{kiendl2020}, a partitioned scheme for immersed shells was recently proposed, in which the relaxation is applied to the beam's acceleration before it is handed to the fluid solver. Motivated by the need to treat the nonlinearity contained in the interaction force applied to the beam, we relax the interaction force instead, and apply the relaxed force $\fbiforceveck = \omega\bfbiforceveck + \left(1-\omega\right)\bfbiforceveckold$ in the proper place, as visualized in Figure~\ref{fig:algo}, instead. See also~\cite{kuttler2006} for another application of a variant of the Aitken method based on force relaxation, and a detailed description of the computation of the relaxation parameter $\omega$.

The computation of the relaxation parameter $\omega$ via the Aitken method defaults to a simple vector-vector product,
making the method highly efficient with regard to the computational cost of the relaxation parameter.
However, even though the Aitken relaxation method works well for many FSI problems, convergence is in general not guaranteed~\cite{kuttler2008fixed}.

\subsubsection{Quasi-Newton-Krylov Solver}
\label{subsubsec:kryloc}

\begin{figure}
\centering

\begin{tikzpicture}[font=\scriptsize]

\tikzstyle{myline} = [-latex,thick]

\draw [fill=gray!20] (-0.2,7.5) node [below right] {Time loop} rectangle (14.2,-1.8);

\draw [fill=gray!60] (0.2,7.0) node [below right] {Quasi-Newton loop} rectangle (14.0,0.2);

\draw [fill=gray!90] (4.5,6.5) node [below right] {Krylov iteration} rectangle (13.8,3.7);

\draw [thick,fill=white] (0.8,6.0) rectangle (4.2,5.0) node [pos=0.5] {Compute $\jac$};

\draw [myline] (4.2,5.5) -- (7.0,5.5);
\draw [thick,fill=white] (7.0,6.0) rectangle (13.6,5.0) node [pos=0.5] {$\alpha^m = \min\limits_{\Delta \bfbiforcevecnewton{k}{m} \in \mathcal{K}^{fbi}_m} \left\|\jac \Delta\bfbiforcevecnewton{k}{m} + \fbirhsnewton{k}{m} \right\|$};

\draw [myline] (11.8,5.0) -- (11.8,4.6);
\draw [thick,fill=white] (10.8,4.1) -- (11.8,4.6) -- (12.8,4.1) -- (11.8,3.6) -- cycle;
\node at (11.8,4.1) {\scriptsize$\alpha_m < tol_{kryl}$};

\draw [myline] (10.8,4.1) node [below] {no} -- (8.7,4.1) node [above,pos=0.5] {$m\leftarrow m+1$};
\draw [ultra thick,fill=colUniBwOr!30] (5.3,4.4) rectangle (8.7,3.8) node [pos=0.5] {Compute $\fbirhsnewton{k}{m}$};
\draw [myline] (5.3,4.1) -- (5.0,4.1) -- (5.0,5.5);

\draw [myline] (11.8,3.6) node [below right] {yes} -- (11.8,3.1) node [left,pos=0.5] {$k\leftarrow k+1$};
\draw [thick,fill=white] (10.0,3.1) rectangle (13.6,2.0) node [pos=0.5] {\begin{tabular}{c}
$\bfbiforcevecnewton{k}{0} = \bfbiforcevecnewton{k-1}{m}$\\
$\fbirhsnewton{k}{0} =  \fbirhsnewton{k-1}{m}$
\end{tabular}};
\draw [myline] (11.8,2.0) -- (11.8,1.5);

\draw [thick,fill=white] (10.2,0.7) -- (11.8,1.5) -- (13.4,0.7) -- (11.8,-0.1) -- cycle;
\node at (11.8,0.7) {\scriptsize$\left\|\fbirhsnewton{k}{0}\right\|<tol_{qn}$};

\draw [myline] (10.2,0.7) node [below] {no} -- (0.5,0.7) -- (0.5,5.5);

\draw [myline] (11.8,-0.1) node [right=6pt] {yes} -- (11.8,-0.5);
\draw [thick,fill=white] (10.8,-1.0) -- (11.8,-0.5) -- (12.8,-1.0) -- (11.8,-1.5) -- cycle;
\node at (11.8,-1.0) {\scriptsize$t < T$};

\draw [myline] (10.8,-1.0) node [below=2pt] {yes} -- (7.6,-1.0) node [pos=0.5,below] {$t\leftarrow t+\Delta t$};
\draw [thick,fill=white] (4.0,-1.5) rectangle (7.6,-0.5) node [pos=0.5] {$\bfbiforcevecnewton{0}{0} = \bfbiforcevecnewton{k}{m}$};
\draw [myline] (4.0,-1.0) -- (0.0,-1.0) node [pos=0.5,below] {$k\leftarrow0, m\leftarrow0$} -- (0.0,5.5);
\draw [myline] (11.8,-1.5) node [right=4pt] {no} -- (11.8,-2.3);

\draw [thick,fill=gray!20] (11.3,-2.3) rectangle (12.3,-2.9) node [pos=0.5] {End};

\draw [thick,fill=gray!20] (-1.6,5.2) rectangle (-0.6,5.8) node [pos=0.5] {Start};
\draw [myline] (-0.6,5.5) -- (0.8,5.5);

\end{tikzpicture}
\caption{Quasi-Newton Method: The computation of~$\fbirhsnewton{k}{m}$ (depicted in the bold-bordered orange box) executes the step ''coupling iteration in each time step'' from the partitioned algorithm depicted in Figure~\ref{fig:algo}.}
\label{fig:quasinewtonalgo}
\end{figure}
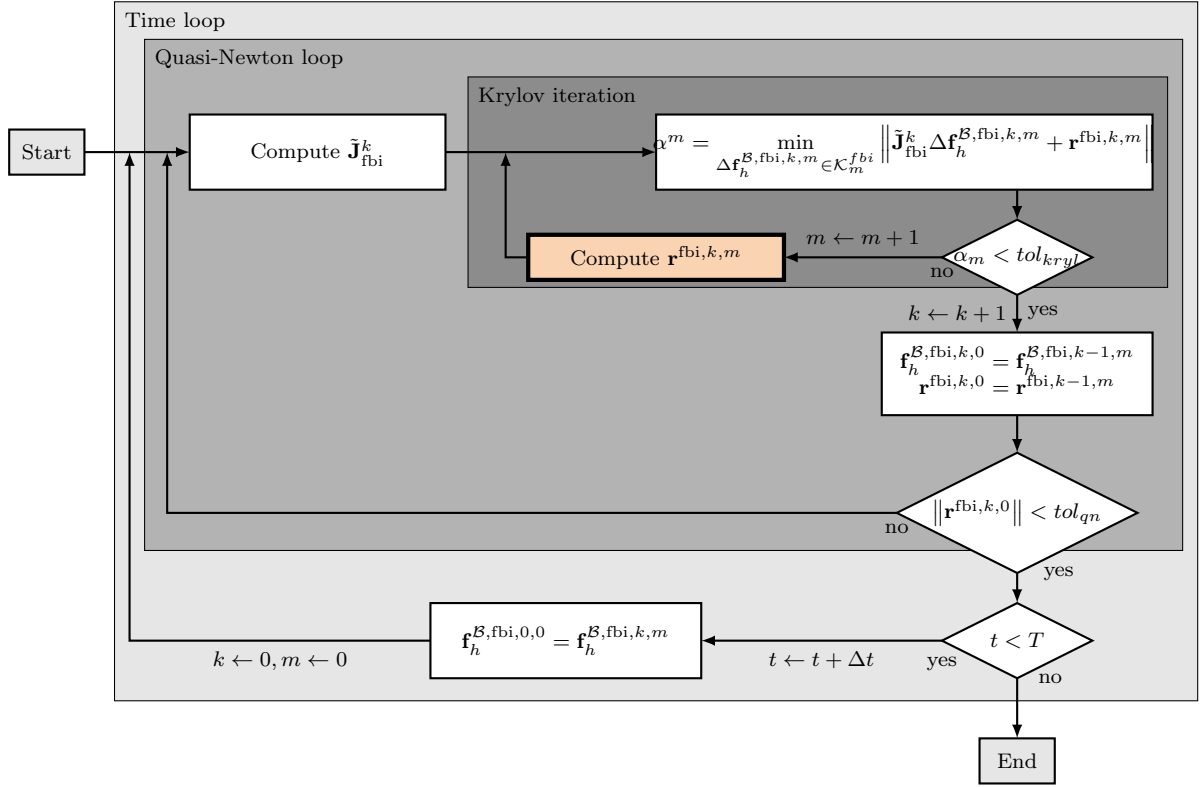

Motivated by the fact that proper treatment of the interface nonlinearity may reduce the number of coupling iterations, we choose a Quasi-Newton scheme applied to the residual
\begin{equation}
	\label{eq:fbi_residual}
	\fbirhsk:=\bfbiforceveck-\bfbiforceveckold.
\end{equation}
This means that we will solve the equation
\begin{equation}
	\label{eq:fbi_newton}
	\jac\Delta \bfbiforceveck  = - \fbirhsk,
\end{equation}
where $\jac$ is an approximation of the linearization of~\eqref{eq:fbi_residual}. To solve the linear system in~\eqref{eq:fbi_newton}, a direct or an iterative linear solver can be used. We choose an iterative Krylov solver, as in that case only the effect of the approximated Jacobian on the current residual, i.e. a matrix-vector product, needs to be known instead of the overall matrix. Based on~\cite{kuttler2008fixed}, we can approximate the effect of the sensitivity on a vector $\mathbf{y}$ as
\begin{equation}
	\label{eq:fe_jacobian}
	\jac\mathbf{y}\approx\dfrac{\mathcal{F}\left(\mathcal{B}\left(\bfbiforceveck+\gamma \left(\gamma +\dfrac{\left|\bfbiforceveck\right|}{\fbirhsk}\right)\mathbf{y}\right)\right)-\bfbiforceveck-\gamma \left(\gamma +\dfrac{\left|\bfbiforceveck\right|}{\fbirhsk}\right)\mathbf{y}-\fbirhsk}{\gamma \left(\gamma +\dfrac{\left|\bfbiforceveck\right|}{\fbirhsk}\right)}.
\end{equation}
Here, $\gamma$ is a user-selected parameter determining the finite differencing step size. The effect of the choice of $\gamma$ will be discussed in detail in Section~\ref{subsec:example_solver}.

This approximate Jacobian can then be used to find the best solution within the Krylov subspace
\begin{equation*}
	\mathcal{K}^{\indFBI}_m := \bfbiforcevectemplate{0}+\mathcal{K}_m\left(\jac, \fbirhstemplate{0}\right).
\end{equation*}
Compared to a direct solver, this choice of an iterative solver for~\eqref{eq:fe_jacobian} comes with the drawback that the residual~\eqref{eq:fbi_residual} will have to be built $m$ times per Quasi-Newton step, where $m$ is the number of iterations of the Krylov solver, as visualized in Figure~\ref{fig:quasinewtonalgo}. In each step of the Krylov solver, the minimization problem
\begin{equation*}
\alpha^m := \min\limits_{\Delta \bfbiforcevecnewton{k}{m} \in \mathcal{K}^{\indFBI}_m} \left\|\jac \Delta\bfbiforcevecnewton{k}{m} + \fbirhsnewton{k}{m} \right\|,
\end{equation*}
is solved until the Krylov subspace $\mathcal{K}^{\indFBI}_m$ is large enough such that $\alpha^m$ falls under a user-prescribed tolerance.
In the case at hand, the residual evaluation in each Krylov step is quite costly as it involves a full execution of the FBI algorithm presented in Figure~\ref{fig:algo}.

On the other hand, no full matrix has to be built or stored, and~\eqref{eq:fe_jacobian} can be directly used to approximate the effect of the sensitivity on the residuals, leading to a matrix-free (Quasi-)Newton-Krylov solver.
For a comprehensive introduction to Krylov subspace methods and more information on linear solvers, the interested reader is referred to the literature, e.g.~\cite{braess2007, saad2003}.

\subsection{Computational aspects}
\label{sec:NonlinearChallenges}

From the definition of the coupling matrices~$\mathbf{D}$ and~$\mathbf{M}$ in~\eqref{eq:mortarmatrices},
it becomes evident that~$\mathbf{M}$ depends on the projection~$\chi$ between the beam and fluid meshes.
Due to the relative motion of the beams {\wrt} to the background fluid mesh,
this projection and consequently the coupling matrix~$\mathbf{M}$ need to be re-computed in every coupling iteration (see also Remark~\ref{remark:segmentation}).

The projection~$\chi$ depends solely on the intersection of the beam mesh with the background fluid mesh.
For its evaluation, a geometrical search to identify pairs of beam and fluid elements interacting with each other
as well as a successive segmentation to generate integration cells for the evaluation of~\eqref{eq:mortarmatrices} has to be performed within every FBI iteration.

In particular for large numbers of immersed beams, and in the context of parallel computing,
this gives rise to the need for computationally efficient parallel search algorithms to find the current beam position within each step of the algorithm.
In this work, the parallel search is performed using a geometrically motivated binning-based communication between processes~\cite{Plimpton1995a}, that follow a distributed memory paradigm,
and a subsequent octree-based search on each shared-memory unit.
The binning strategy represents a geometric partitioning method that is tailored to models where interactions between high numbers of bodies have to be resolved
as it is commonly also the case for the simulation of particle methods~\cite{rhoades1992,tang2014,eichinger2021}.
We refer the interested reader to~\cite{eichinger2021,Mayr2023a} for a detailed description of the employed binning implementation and a discussion of its advantages and special features.

\section{Numerical Examples}
\label{sec:examples}
The following numerical examples are chosen to investigate the behavior of the introduced numerical and algorithmic building blocks and to demonstrate the performance of the proposed methodology also for complex problems with large interface displacements and many interacting fibers. For all examples, the fluid is assumed to be initially at rest. All models were created using the pre-processor MeshPy~\cite{MeshPyWebsite}, and all simulations were executed using the multi-physics research code BACI~\cite{baci}.

\subsection{Single elastic beam}
\label{subsec:example_oscilating}

\begin{figure}[b]
	\begin{center}
		\begin{subfigure}{0.6\textwidth}
			\resizebox{\textwidth}{!}{




\begin{tikzpicture}

\tikzstyle{measline} = [latex-latex]
\tikzstyle{axisline} = [-latex]

\tikzstyle{beam} = [ultra thick, gray]

\fill[white,opacity=.5] (0,0,0)-- (3,0,0) -- (3,1,0)  -- (0,1,0) --cycle;
\draw [thick] (3,0,0)  -- (3,1,0) --(0,1,0) -- (0,0,0) -- (3, 0, 0);

\draw (0,-0.05,0) -- (0,-0.25,0);
\draw (3,-0.05,0) -- (3,-0.25,0);
\draw[measline] (0,-0.2,0)-- (3,-0.2,0) node[pos=.5, below] {\tiny$l=3$};

\draw (-0.05,0,0) -- (-0.25,0,0);
\draw (-0.05,1,0) -- (-0.25,1,0);
\draw[measline] (-0.2,0,0)-- (-0.2,1,0) node[pos=.5, left] {\tiny$h=1$};

\draw [beam] (1.5,0,0) -- (1.5, 0.5 ,0);
\draw (1.55,0.5,0) -- (1.75, 0.5 ,0);
\draw[measline] (1.7,0, 0) -- (1.7, 0.5 ,0) node[pos=.5,right] {\tiny$h_\indBeam=0.5$};

\node at (0.2,0.5) {\tiny{$\Gamma_\mathrm{in}$}};
\node at (3.3,0.5) {\tiny{$\Gamma_\mathrm{out}$}};

\draw[axisline] (-0.5, -0.5, 0) -- (-0.1, -0.5, 0) node[below] {\tiny$x$};
\draw[axisline] (-0.5, -0.5, 0) -- (-0.5, -0.1, 0) node[left] {\tiny$y$};

\end{tikzpicture}

}
			\subcaption{Geometric setup of the quasi-2D example of an immersed elastic beam}
			\label{fig:biegebalken:geometry}
		\end{subfigure}
		\begin{subfigure}{0.38\textwidth}
			\resizebox{\textwidth}{!}{\input{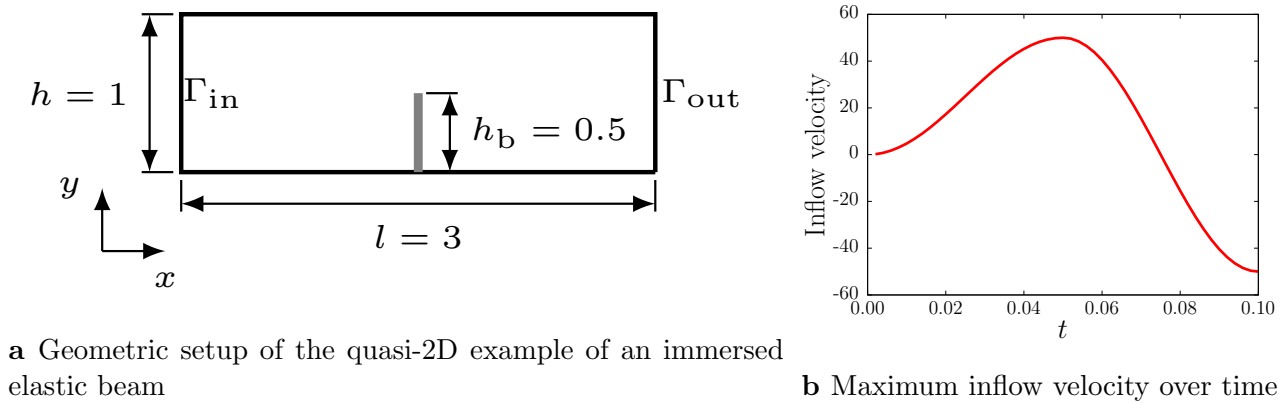}}
			\subcaption{Maximum inflow velocity over time}
			\label{fig:biegebalken_inflow_vel}
		\end{subfigure}
	\end{center}
	\caption{Problem setup for the single elastic beam immersed in a fluid channel}
\end{figure}
\begin{figure}
	\begin{subfigure}{0.455\textwidth}
		\includegraphics[width=\textwidth]{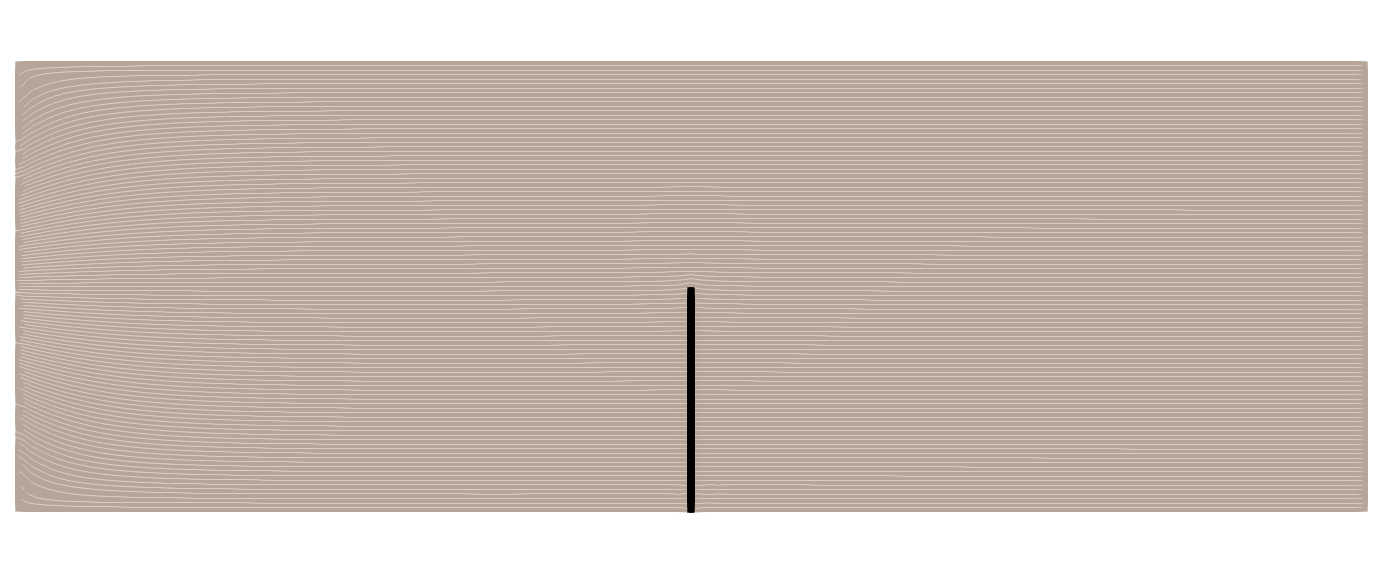}
		\subcaption{Initial configuration}
		\label{fig:biegebalken_initial}
	\end{subfigure}
	\begin{subfigure}{0.535\textwidth}
		\includegraphics[width=\textwidth]{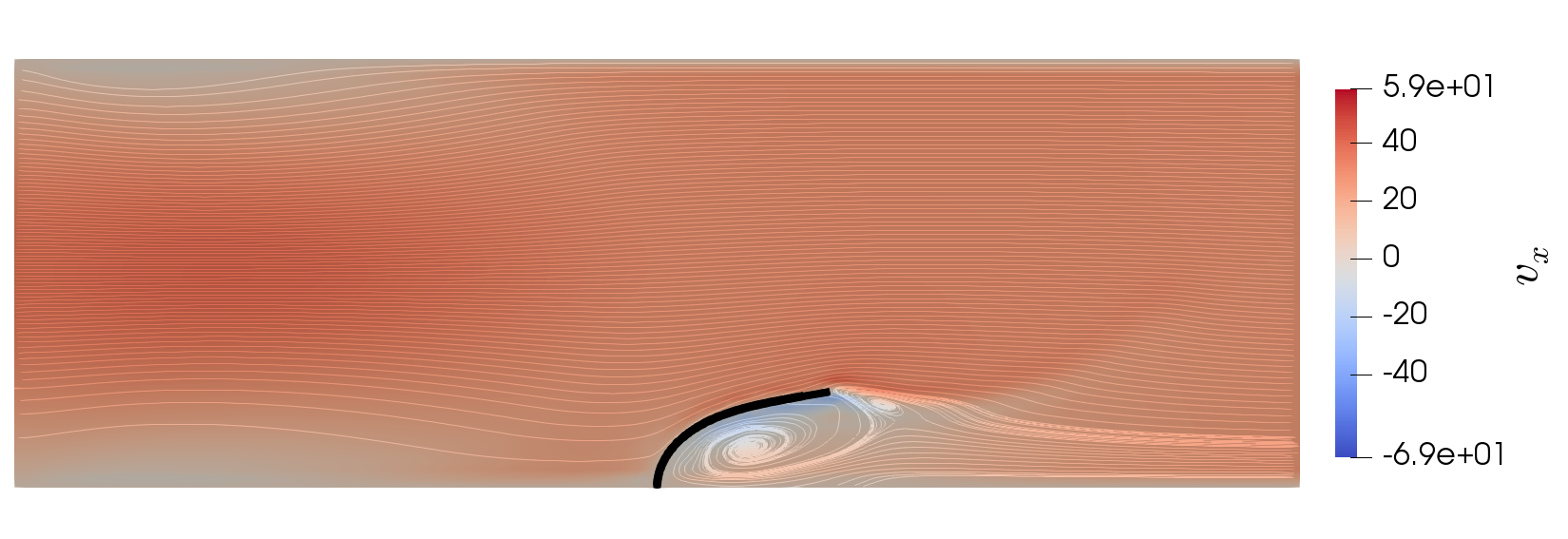}
		\subcaption{Time $t=0.062$}
		\label{fig:biegebalken_b}
	\end{subfigure}
	\begin{subfigure}{0.455\textwidth}
		\includegraphics[width=\textwidth]{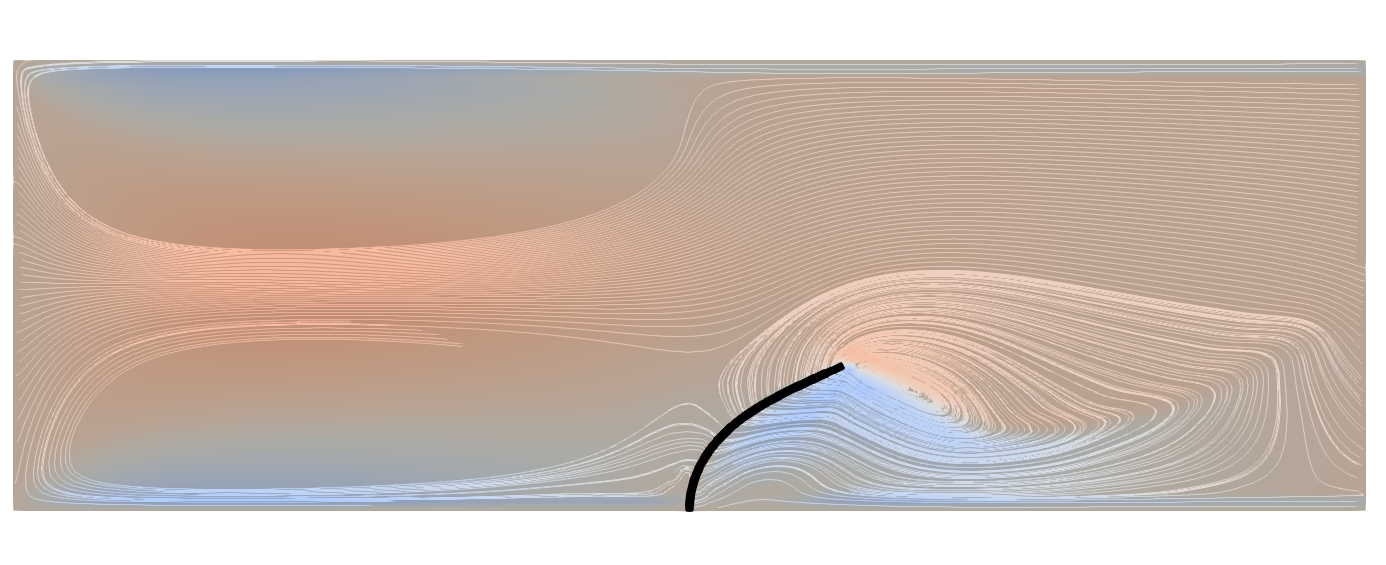}
		\subcaption{Time $t=0.074$}
		\label{fig:biegebalken_c}
	\end{subfigure}
	\begin{subfigure}{0.535\textwidth}
		\includegraphics[width=\textwidth]{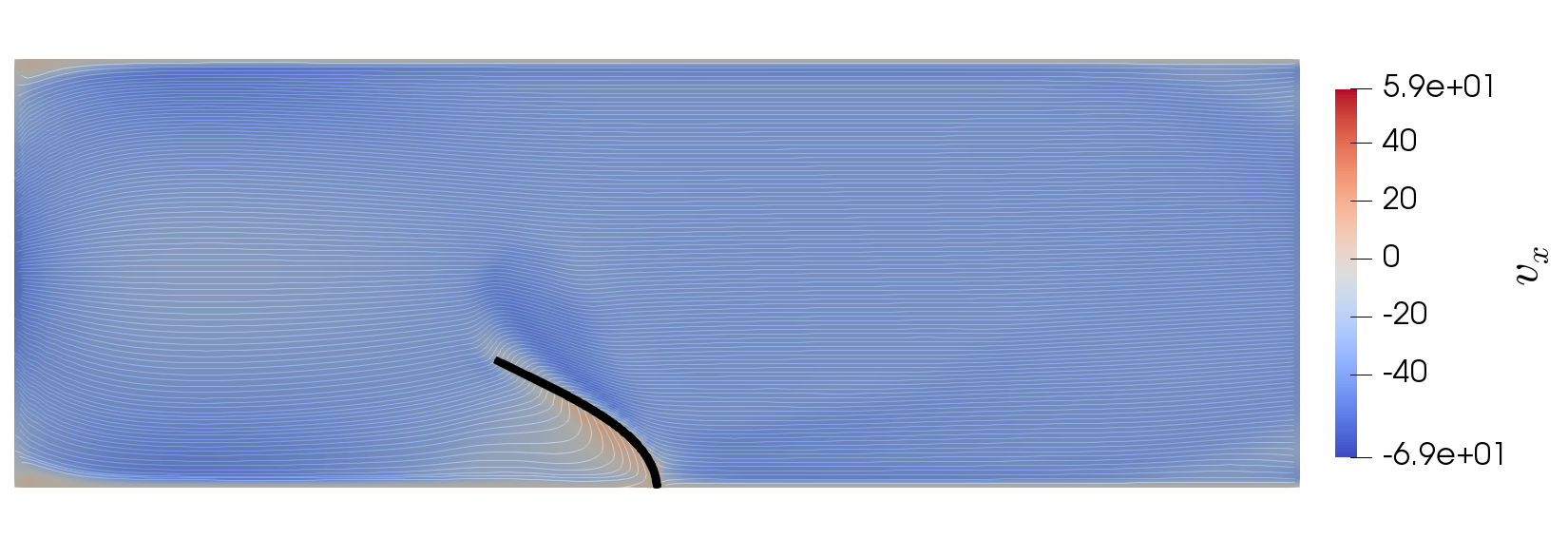}
		\subcaption{Time $t=0.1$}
		\label{fig:biegebalken_d}
	\end{subfigure}
	\caption{Velocity solution in channel direction at different time steps}
	\label{fig:biegebalken_picture}
\end{figure}

The intention of this example is to visualize the effect of the method's numerical parameters, namely penalty parameter, constraint enforcement technique, and, in the case of the mortar-type approach, Lagrange multiplier shape functions, on the solution of a fluid-beam system exhibiting large beam displacements. For this purpose, the problem is set up as pseudo quasi-2D in order to facilitate large displacements. The geometric setup of the channel and the beam is given in Figure~\ref{fig:biegebalken:geometry}. The fluid inflow is prescribed on the left end of the fluid channel as a parabolic flow profile with respect to the channel height, and oscillating in time. The time evolution of the fluid velocity at the middle of the channel height is visualized in Figure~\ref{fig:biegebalken_inflow_vel}. inspired by the real physical properties of water, the fluid has the density $\rho^\indFluid = 1$, and the viscosity $\nu^\indFluid=0.004$. The channel has a height of $1$, the length is $3$, and the depth is $0.06$ over two fluid elements. To allow for negative fluid velocities also on the Neumann boundary on the right, backflow boundary conditions, as analyzed in~\cite{bertoglio2018}, are used. Non-penetration conditions are applied to all other channel surfaces. The immersed beam has a height $h_\indBeam=0.5$, cross-sectional area $A=10^{-4}$, density $\rho^\indBeam=10$, Young's modulus $E^{\indBeam}=10^7$, and is modeled using a hyperelastic material. Figure~\ref{fig:biegebalken_picture} illustrated the geometrical configuration as well as the fluid velocity in channel direction at different time snippets for the mortar-type method with linear Lagrange multiplier shape functions and a penalty parameter $\epsilon=10^4$. The initial configuration can be found in Figure~\ref{fig:biegebalken_initial}. Figure~\ref{fig:biegebalken_b} demonstrates the beam's deflection to the right, just before the flow direction changes, and  Figures \ref{fig:biegebalken_c}- \ref{fig:biegebalken_d} show the beam's deflection after a change in flow direction.

\begin{figure}
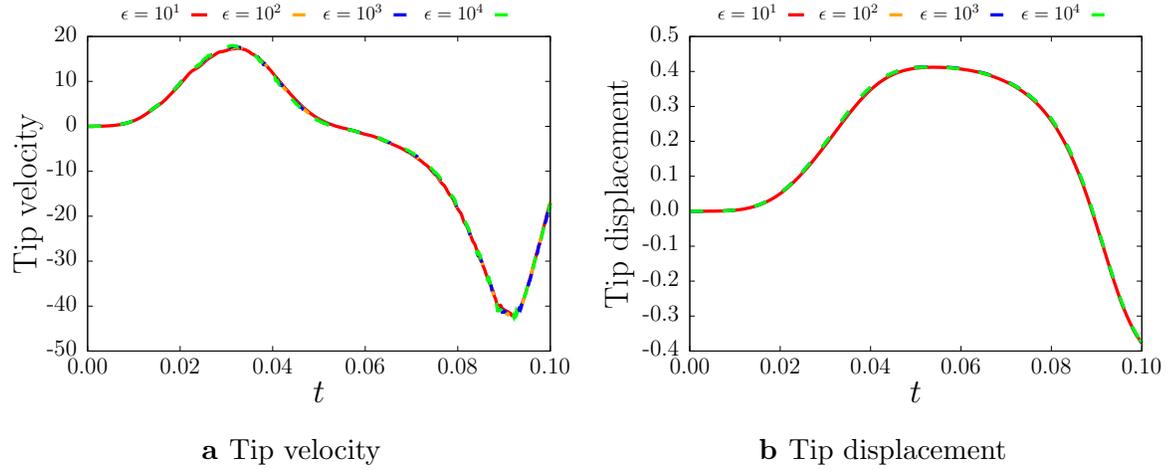

	\begin{subfigure}{0.45\textwidth}
		\resizebox{\textwidth}{!}{\input{pictures/twoway_oscillating/vel_over_penalty_line2.tex}}
		\subcaption{Tip velocity}
	\end{subfigure}
	\begin{subfigure}{0.45\textwidth}
		\resizebox{\textwidth}{!}{\input{pictures/twoway_oscillating/disp_over_penalty_line2.tex}}
		\subcaption{Tip displacement}
	\end{subfigure}
	\caption{Velocity and displacement for the mortar-type approach with linear Lagrange multiplier shape functions and different penalty parameters $\epsilon$}
	\label{fig:biegebalken_mortar}
\end{figure}

To analyze the effect of the penalty parameter on the simulation result, Figure~\ref{fig:biegebalken_mortar} demonstrates the beam tip velocity and tip displacement for a solution obtained with linear Lagrange multiplier shape functions for various penalty parameter values.
Figure~\ref{fig:biegebalken_mortar} demonstrates that, even for relatively small penalty parameters, the tip velocity as well as tip displacement is adequately captured without large mismatches between the different solutions over the whole duration of the simulation. The maximum relative differences with respect to time for the tip displacement $d^\text{tip}_{10^4}$, computed using $\epsilon=10^4$ compared to the tip displacement computed with lower penalty parameters, are shown in Table \ref{tab:penalty_errors}. As expected, for growing penalty values, this difference becomes smaller suggesting convergence of the solution with respect to the penalty parameter. Furthermore, for the examined penalty values, the maximum difference stays below 2.5\% even for the relatively large displacements and high velocities exhibited by the model problem.

Similarly, Table \ref{tab:lagrangediscret_errors} contains the comparison of the tip displacement computed using $\epsilon=10^2$ and cubic Lagrange multiplier shape functions with the solution for linear and quadratic shape functions. For the analyzed example, the difference introduced by the Lagrange shape functions stays well below 1\%.

\begin{table}
	\begin{center}
		\begin{tabular}{|c|c|c|}
			$\dfrac{\left|\left|d^\text{tip}_{10^4}\left(t\right)-d^\text{tip}_{10^1}\left(t\right)\right|\right|_{\infty}}{\left|\left|d^\text{tip}_{10^4}\left(t\right)\right|\right|_{\infty}}$ & $\dfrac{\left|\left|d^\text{tip}_{10^4}\left(t\right)-d^\text{tip}_{10^2}\left(t\right)\right|\right|_{\infty}}{\left|\left|d^\text{tip}_{10^4}\left(t\right)\right|\right|_{\infty}}$ & $\dfrac{\left|\left|d^\text{tip}_{10^4}\left(t\right)-d^\text{tip}_{10^3}\left(t\right)\right|\right|_{\infty}}{\left|\left|d^\text{tip}_{10^4}\left(t\right)\right|\right|_{\infty}}$ \\
			\\
			\hline \\
			2.33\% & 1.23\% & 0.49\%
		\end{tabular}
		\caption{Relative differences introduced by the penalty parameter}
		\label{tab:penalty_errors}
	\end{center}
\end{table}

\begin{table}
	\begin{center}
		\begin{tabular}{|c|c|}
			$\dfrac{\left|\left|d^\text{tip}_\text{cub}\left(t\right)-d^\text{tip}_\text{quad}\left(t\right)\right|\right|_{\infty}}{\left|\left|d^\text{tip}_\text{cub}\left(t\right)\right|\right|_{\infty}}$ & $\dfrac{\left|\left|d^\text{tip}_\text{cub}\left(t\right)-d^\text{tip}_\text{lin}\left(t\right)\right|\right|_{\infty}}{\left|\left|d^\text{tip}_\text{cub}\left(t\right)\right|\right|_{\infty}}$ \\
			\\
			\hline \\
			0.48\% & 0.1\%
		\end{tabular}
		\caption{Relative differences introduced by the Lagrange multiplier shape functions}
		\label{tab:lagrangediscret_errors}
	\end{center}
\end{table}
Even though not necessarily negligible, the above reported numerical effects are well within the range of the modeling error reported in~\cite{zunino2016}. There, a beam based stent model was compared with a continuum mechanics based reference model, and a maximum displacement error of 4\% was found. The modeling error for the proposed FBI approach will be further investigated in Section \ref{subsec:example_beam_comparison}. 

In conclusion, the represented example exploring the effect of numerical parameter choices on the simulation results shows convergence of the solution with respect to the value of the penalty parameter as well as the Lagrange multiplier shape functions. Besides convergence, the results suggest that also the use of moderate penalty parameters as well as linear Lagrange multiplier shape functions allows to sufficiently capture the overall solution even for large displacements. That constitutes a desirable result as these simple choices ease the solution procedure of the fields' resulting linear systems of equations due to the avoidance of ill-conditioning effects.

\subsection{Comparison to 3D reference solution}
\label{subsec:example_beam_comparison}

\begin{figure}
	\begin{center}
		\resizebox{0.9\textwidth}{!}{




\begin{tikzpicture}

\tikzstyle{measline} = [latex-latex]
\tikzstyle{axisline} = [-latex]

\fill[white,opacity=.5] (0,0,0)-- (3,0,0) -- (3,1,0)  -- (0,1,0) --cycle;
\fill[white,opacity=.5] (0,0,1)-- (3,0,1) -- (3,1,1)  -- (0,1,1) --cycle;
\fill[white,opacity=.5] (0,1,0)-- (0,1,1) -- (3,1,1) -- (3,1,0)--cycle;
\fill[white,opacity=.5] (0,0,0)-- (0,0,1) -- (3,0,1) -- (3,0,0)--cycle; 
\draw[] (0,0,1) -- (3,0,1) -- (3,1,1) --(0,1,1) --(0,0,1)
(3,0,1) -- (3,0,0)  -- (3,1,0) --(0,1,0) -- (0,1,1)    
(3,1,1) -- (3,1,0);
\draw[dashed] (0,0,0) -- (0,0,1) (0,0,0)-- (3,0,0) (0,0,0)-- (0,1,0);

\draw[measline] (0.2,0,1.5)-- (3.2,0,1.5) node[pos=.5, below] {\tiny$l=3$};
\draw[measline] (3.2,0,0.1)-- (3.2,0,1.1) node[pos=.7, right] {\tiny$b=1$};
\draw[measline] (-0.15,0,1)-- (-0.15,1,1) node[pos=.5, left] {\tiny$h=1$};

\node (A) [cylinder,draw=black,aspect=0.4,
minimum height=0.6cm,minimum width=0.000001cm,
shape border rotate=90,fill=gray!50] at (1.5, 0.2, 0.5 ) {};

\draw[dashed]
let \p1 = ($ (A.after bottom) - (A.before bottom) $),
\n1 = {0.5*veclen(\x1,\y1)},
\p2 = ($ (A.bottom) - (A.after bottom)!.5!(A.before bottom) $),
\n2 = {veclen(\x2,\y2)}
in
(A.before bottom) arc [start angle=0, end angle=180,
x radius=\n1, y radius=\n2];

\draw (1.5,0, 0.5) -- (1.5, 0.5 ,0.5);

\draw[measline] (1.725,0, 0.5) -- (1.725, 0.5 ,0.5) node[pos=.7, anchor = west] {\hspace{-0.05cm}\tiny$h_\indBeam=0.5$};

\node at (-0.2,0.3) {\tiny{$\Gamma_\mathrm{in}$}};
\node at (3.3,0.3) {\tiny{$\Gamma_\mathrm{out}$}};

\draw[axisline] (-2.5, -1, -1.3) -- (-2.0, -1,-1.3) node[below] {\tiny$x$};
\draw[axisline] (-2.5, -1, -1.3) -- (-2.5, -0.5, -1.3) node[left] {\tiny$y$};
\draw[axisline] (-2.5, -1, -1.3) -- (-2.5, -1, -0.6) node[left] {\tiny$z$};

\end{tikzpicture}

}
		\caption{Elastic obstacle immersed in a fluid channel}
		\label{fig:obstaclesetup}
	\end{center}
\end{figure}
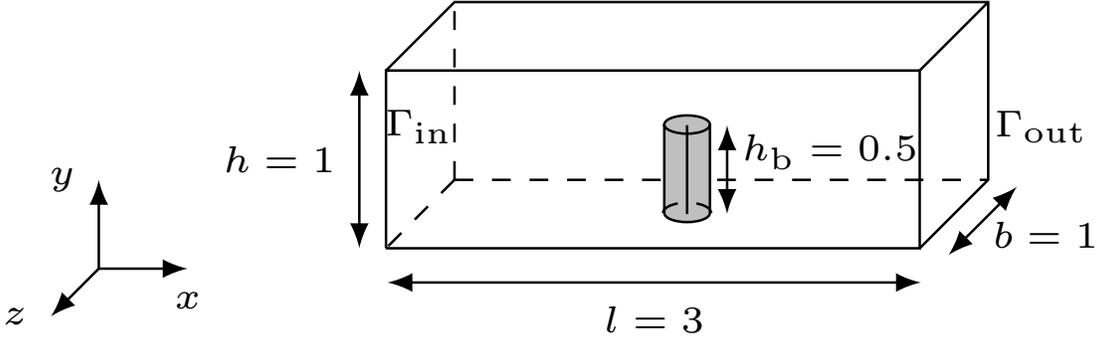

Within this section, the spatial convergence behavior of the proposed fully-coupled mixed-dimensional method compared to a fully resolved ALE-based FSI approach as presented in~\cite{Kloppel2012, mayr2014} is analyzed. To this end, a fluid channel with the dimensions $3\times 1\times 1$ is examined. An immersed beam of height $0.5$, and with radius $0.01$, is fixed in the middle of the bottom of the channel as sketched in Figure~\ref{fig:obstaclesetup}. The beam material is modeled using a hyperelastic material with the Young's modulus $E^{\indBeam} = E^{\mathrm{s}}=5\cdot 10^9$ and the density $\rho^{\indBeam}=\rho^{\mathrm{s}}=10^1$ for both, the beam as well as the fully resolved model. Meanwhile, the fluid is modeled as a Newtonian fluid with different densities $\rho^\indFluid$ and various constant viscosities $\nu^\indFluid$. The simulations are run using the Generalized-$\alpha$ time integration method with a time step size $\Delta t=10^{-4}$ and a spectral radius $\rho_\infty=1$ for the structure fields~\cite{genalpha}, and a One-step-$\theta$ method with $\theta=1.0$ for the fluid field. Within the first 500 steps, the inflow velocity is slowly accelerated until a parabolic inflow profile
\begin{equation*}
	\mathbf{v}_{\mathrm{in}} = 1600\cdot y \cdot \left(1-y\right)\cdot z \cdot \left(1-z\right)
\end{equation*}
is reached.
A zero traction boundary condition is prescribed on the outflow boundary and no-slip boundary conditions are enforced on all other sides.
The FBI simulations are run with a penalty parameter $\epsilon=10^4$ and linear Lagrange multiplier shape functions. The fully resolved 3D reference solution, in turn, is simulated with the Lagrange multiplier based ALE approach presented in~\cite{mayr2014, Kloppel2012}, using 2,183,072 and 10,240 linear hexahedral finite elements for the fluid and beam fields, respectively.

\begin{figure}[h!]
	\begin{subfigure}[b]{0.48\textwidth}
		\includegraphics[width=\textwidth]{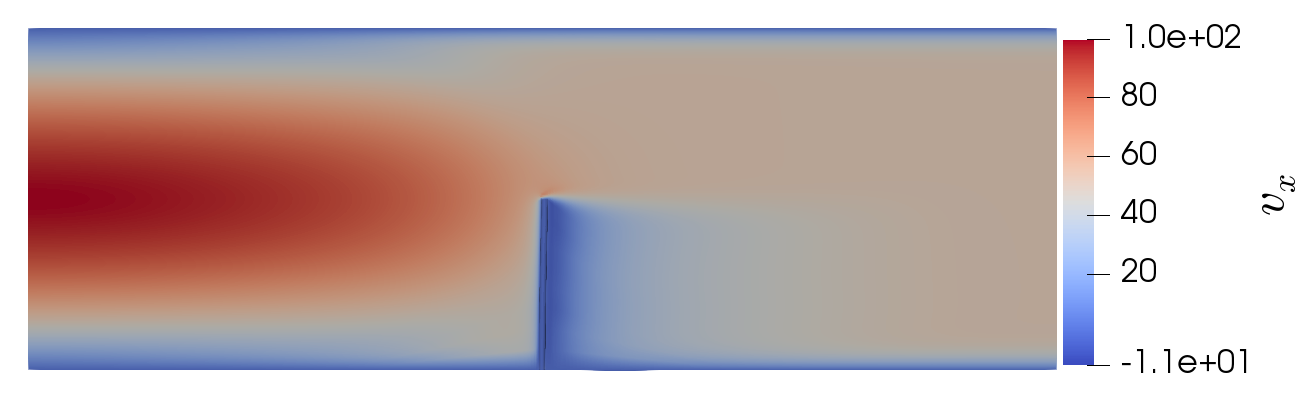}
		\vspace{0.6cm}
		\includegraphics[width=\textwidth]{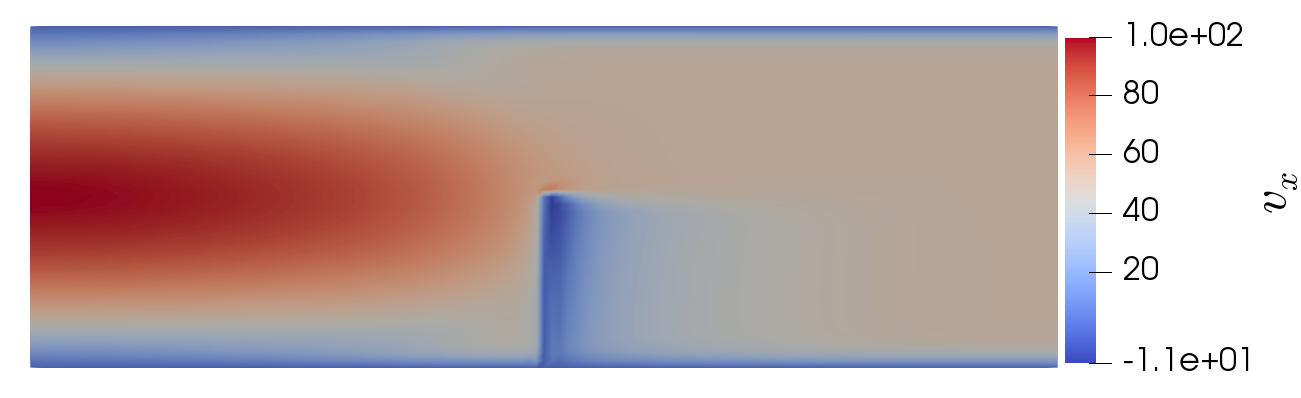}
		\subcaption{Velocity of the reference solution (top), and the FBI solution with the fluid element size equal to the beam diameter}
		\label{subfig:comparison}
	\end{subfigure}
	\hspace{0.3cm}
	\begin{subfigure}[b]{0.48\textwidth}
		\resizebox{\textwidth}{!}{\input{pictures/biegebalken_comparison/new/convergencev4.tex}}
		\subcaption{Behavior of the relative L2 error of the fluid velocity under uniform mesh refinement for varying Reyolds number}
		\label{subfig:convergence}
	\end{subfigure}
	\caption{Convergence behavior under uniform mesh refinement}
\end{figure}

The obtained velocity solution for the FBI example with a fluid element size $\dfrac{1}{44}$ and the solution of the reference simulation with $\rho^{\indFluid}=1.0$ and $\nu^{\indFluid}=0.032$ are depicted in Figure~\ref{subfig:comparison}.

Figure~\ref{subfig:convergence} visualizes the relative L2 error of the fluid velocity solution compared to the fully resolved 3D reference solution under uniform mesh refinement and for different fluid material configurations at time $t=0.5$. Linear convergence can be observed for coarser meshes, which is in agreement with the results in~\cite{hagmeyer2022}, where linear convergence was recovered for a fixed beam obstacle.

In addition to the spatial discretization error, the shown L2 error contains the modeling difference introduced by using beam theory instead of continuum mechanics. Thus, spatial convergence can only be expected, firstly, as long a the modeling assumptions regarding the beam diameter are fulfilled, and secondly, as long as the error due to the spatial fluid field discretization is larger than the aforementioned modeling error.

While the convergence under uniform mesh refinement is limited by the modeling assumptions, it is important to note that the solution obtained with the presented mixed-dimensional approach converges against the computed 3D reference solution down to a point where the fluid element size nearly matches the beam diameter. In particular, no qualitative dependency of the convergence rate on the fluid material parameters can be observed. The L2-error for the optimum of any of the chosen material parameter sets falls under 1.5\% while, at the same time, the mixed-dimensional approach reduces the number of DoFs to a fraction of $\dfrac{1}{25}$, {\ie} by 96\%.

While a fully resolved 3D surface-coupled model is still recommended for applications where phenomena near the interface are of interest, this example shows that the solution obtained by the mixed-dimensional method converges against the reference solution and provides a globally matching solution at a fraction of the computational cost, thus validating the proposed methodology in a quite remarkable manner.

\subsection{Comparison of partitioned solver strategies}
\label{subsec:example_solver}

\begin{figure}[h!]
	\begin{center}
		\begin{subfigure}{0.45\textwidth}
			\resizebox{\textwidth}{!}{\input{pictures/solver_comparison/biegebalken/lambda_ite1_totaltime.tex}}
			\subcaption{Using linearization to compute the finite differences approximation}
			\label{fig:solver_parameter_stiff_linear_time}
		\end{subfigure}
		\begin{subfigure}{0.45\textwidth}
			\resizebox{\textwidth}{!}{\input{pictures/solver_comparison/biegebalken/lambda_ite5_totaltime.tex}}
			\subcaption{Solving the full nonlinear problems to compute the finite differences approximation}
			\label{fig:solver_parameter_stiff_nonlinear_time}
		\end{subfigure}
		\begin{subfigure}{0.45\textwidth}
			\resizebox{\textwidth}{!}{\input{pictures/solver_comparison/biegebalken/lambda_ite1_iterations.tex}}
			\subcaption{Using linearization to compute the finite differences approximation}
			\label{fig:solver_parameter_stiff_linear_iteration}
		\end{subfigure}
		\begin{subfigure}{0.45\textwidth}
			\resizebox{\textwidth}{!}{\input{pictures/solver_comparison/biegebalken/lambda_ite1_liniterations.tex}}
			\subcaption{Number of FSI evaluations when using the linearized single field systems}
			\label{fig:solver_parameter_stiff_linear_liniteration}
		\end{subfigure}
		\caption{Analysis of the acceleration techniques}
		\label{fig:solver_parameter_stiff}
	\end{center}
\end{figure}

Figures \ref{fig:solver_parameter_stiff_linear_time} and \ref{fig:solver_parameter_stiff_nonlinear_time} show the computational time required for the Aitken relaxation method compared to the MFNK approach with different values of $\gamma$. For the values visualized in Figure~\ref{fig:solver_parameter_stiff_linear_time}, only the linearized single field equations were solved to compute the finite differences approximation of the Jacobian, while the time measurements for the case where the full nonlinear field equations are solved are depicted in \ref{fig:solver_parameter_stiff_nonlinear_time}.
It is notable in the current example that both, linearized and fully nonlinear approaches to approximate the finite differences, yield comparable computing times.

Figures \ref{fig:solver_parameter_stiff_linear_time} and \ref{fig:solver_parameter_stiff_nonlinear_time} reveal that the MFNK solvers solve the considered problem faster than the Aitken relaxation method within the first 200 steps, independently of the used $\gamma$. After 200 steps, the amount of required iterations increases for all MFNK solvers from $2$ to $3$ iterations. 
For the MFNK solver, each iteration contains a full Krylov solve and multiple evaluations of the residual~\eqref{eq:fbi_residual}. Therefore, this additional iteration increases the computation time significantly.
Afterwards, the speed of the solver highly depends on the value of $\gamma$. This is not surprising, as the approximation of the Jacobian becomes better as $\gamma$ tends to zero. As shown in Figure~\ref{fig:solver_parameter_stiff_linear_liniteration}, this is not the main reason for the speedup here, as the number of evaluations of the FSI residuum stays the same.
Nevertheless, for each residual evaluation, the force $\mathcal{F}\left(\mathcal{B}\left(\mathbf{f}_h^{\mathcal{B}, \indFBI}+\delta\mathbf{y}\right)\right)$ has to be evaluated. The smaller $\gamma$, the smaller the step length $\delta\mathbf{y}$, which can often be beneficial for the iterative linear solver applied to the fluid partition.

Looking at Figure~\ref{fig:solver_parameter_stiff_linear_iteration}, it becomes evident that the MFNK solvers continuously stay below three iterations. The Aitken relaxation method, in contrast, shows large variations in the number of iterations, particularly within the first 200 time steps, where the maximum iteration count is 85.

\begin{figure}[h!]
	\begin{center}
		\begin{subfigure}{0.45\textwidth}
			\resizebox{\textwidth}{!}{\input{pictures/solver_comparison/biegebalken_light/light_lambda_ite1_totaltime.tex}}
			\subcaption{Using linearization to compute the finite differences approximation}
			\label{fig:solver_parameter_soft_linear_time}
		\end{subfigure}
		\begin{subfigure}{0.45\textwidth}
			\resizebox{\textwidth}{!}{\input{pictures/solver_comparison/biegebalken_light/light_lambda_ite5_totaltime.tex}}
			\subcaption{Solving the full nonlinear single field systems to compute the finite differences approximation}
			\label{fig:solver_parameter_soft_nonlinear_time}
		\end{subfigure}
		\begin{subfigure}{0.45\textwidth}
			\resizebox{\textwidth}{!}{\input{pictures/solver_comparison/biegebalken_light/light_lambda_ite1_liniterations.tex}}
			\subcaption{Number of FSI evaluations when using the linearized single field systems}
			\label{fig:solver_parameter_soft_linear_liniteration}
		\end{subfigure}
		\begin{subfigure}{0.45\textwidth}
			\resizebox{\textwidth}{!}{\input{pictures/solver_comparison/biegebalken_light/light_lambda_ite5_liniterations.tex}}
			\subcaption{Number of FSI evaluations when solving the full nonlinear single field systems}
			\label{fig:solver_parameter_soft_nonlinear_liniteration}
		\end{subfigure}
		\caption{Analysis of the acceleration techniques for a light structure}
		\label{fig:solver_parameter_soft}
	\end{center}
\end{figure}

Figure~\ref{fig:solver_parameter_soft} visualizes the solver behavior for a variation of the problem above, where the density of the beam is decreased by a factor of 2. For this somewhat more challenging setup, the typically required number of iterations, and with that the computational time of the Aitken relaxation technique, increases significantly compared to the example with higher density.
Here, the MFNK method with any of the considered values for $\gamma$ exhibits preferable computational times compared to the Aitken method when using the linearized single field equations to compute the finite differences approximation, as depicted in Figure~\ref{fig:solver_parameter_soft_linear_time}. In contrast, Figure~\ref{fig:solver_parameter_soft_nonlinear_time} indicates that, for $\gamma=10^{-5}$,  solving the full single field equations leads to significantly higher computational times within the last 300 steps of the simulation. In any case, Figures \ref{fig:solver_parameter_soft_linear_liniteration} and \ref{fig:solver_parameter_soft_nonlinear_liniteration} demonstrate that the number of FSI residuum evaluations required to solve the considered problem using the MFNK method is consistently lower than with the Aitken method, independently of the configuration.

As noted in~\cite{kuttler2008fixed}, the Aitken relaxation method is a fairly simple and cheap method that often performs well in accelerating convergence. Nevertheless, choosing good parameters for the MFNK method can still speed up the convergence considerably. It is noteworthy that at least one parameter also has to be chosen for the Aitken relaxation method: the maximum number of allowed iterations until divergence is assumed. As shown in Figures \ref{fig:solver_parameter_soft_linear_liniteration} and \ref{fig:solver_parameter_soft_nonlinear_liniteration}, the number of iterations can vary considerably, making this choice particularly difficult but integral for the robustness of the simulation.

\subsection{Submerged Vegetation}
\label{subsec:vegetation}

\begin{figure}[h!]
\centering
	\begin{subfigure}{0.6\textwidth}
		\resizebox{\textwidth}{!}{




\begin{tikzpicture}

\tikzstyle{measline} = [latex-latex]
\tikzstyle{axisline} = [-latex]
\tikzstyle{plant}=[ultra thick, gray]

\fill[white,opacity=.5] (0,0)-- (2,0) -- (2,1)  -- (0,1) --cycle;

\draw[thick] (2,0) -- (0,0)  -- (0,1) --(2,1);
\draw[dashed] (2,1) -- (2,0);

\draw[thick] (2.5,0) -- (4.5,0)  -- (4.5,1) -- (2.5,1);
\draw[dashed] (2.5,1) -- (2.5,0) node[pos=.5, anchor=east]{\tiny{...}};

\draw (0,-0.05) -- (0,-0.25);
\draw (1,-0.05) -- (1,-0.25);
\draw[measline] (0,-0.2)-- (1,-0.2) node[pos=.5, below] {\tiny$1$};

\draw (3.3,-0.05) -- (3.3,-0.25);
\draw (4.5,-0.05) -- (4.5,-0.25);
\draw[measline] (3.3,-0.2)-- (4.5,-0.2) node[pos=.5, below] {\tiny$3$};

\draw (0,1.05) -- (0,1.25);
\draw (4.5,1.05) -- (4.5,1.25);
\draw[measline] (0,1.2)-- (4.5,1.2) node[pos=.5, anchor=south] {\tiny$l=13$};

\draw (1.2,-0.05) -- (1.2,-0.45);
\draw (1.4,-0.05) -- (1.4,-0.45);
\draw (1.2,-0.4) -- (1.4,-0.4);
\draw [axisline] (1.0,-0.4) -- (1.2,-0.4);
\draw [axisline] (1.7,-0.4) node [below] {\tiny$\frac{10}{216}$} -- (1.4,-0.4);

\draw (-0.05,0) -- (-0.25,0);
\draw (-0.05,1) -- (-0.25,1);
\draw[measline] (-0.2,0)-- (-0.2,1) node[pos=.5, left] {\tiny$h=1$};

\draw [plant] (1,0) -- (1, 0.3);
\draw [plant] (1.2,0) -- (1.2, 0.3);
\draw [plant] (1.4,0) -- (1.4, 0.3);
\draw [plant] (1.6,0) -- (1.6, 0.3);
\draw [plant] (1.8,0) -- (1.8, 0.3);

\draw [plant] (2.7,0) -- (2.7, 0.3);
\draw [plant] (2.9,0) -- (2.9, 0.3);
\draw [plant] (3.1,0) -- (3.1, 0.3);
\draw [plant] (3.3,0) -- (3.3, 0.3);

\draw (2.05,0) -- (2.25,0);
\draw (1.85,0.3) -- (2.25,0.3);
\draw [axisline] (2.2,0.6) -- (2.2,0.3);
\draw [axisline] (2.2,-0.4) node [right] {\tiny$h_b=0.3$} -- (2.2,0.0);
\draw (2.2,0) -- (2.2,0.3);

\draw[axisline] (-0.7, -0.7 ) -- (-0.3, -0.7) node[below] {\tiny$x$};
\draw[axisline] (-0.7, -0.7) -- (-0.7, -0.3) node[left] {\tiny$y$};

\end{tikzpicture}

}
	\end{subfigure}
	\hspace{1cm}
	\begin{subfigure}{0.27\textwidth}
		\resizebox{\textwidth}{!}{




\begin{tikzpicture}

\tikzstyle{measline} = [latex-latex]
\tikzstyle{axisline} = [-latex]
\tikzstyle{plant}=[ultra thick, gray]

\fill[white,opacity=.5] (0,0) -- (2,0) -- (2,2) -- (0,2) -- cycle;
\draw[thick] (2,0) -- (0,0) -- (0,2) -- (2,2) -- (2, 0);

\draw (0,2.05) -- (0,2.25);
\draw (2,2.05) -- (2,2.25);
\draw[measline] (0,2.2) -- (2,2.2) node[pos=.5, anchor=south] {\small$b=1$};


\draw (0.0,-0.05) -- (0.0,-0.25);
\draw (0.3,-0.05) -- (0.3,-0.25);
\draw[measline] (0.0,-0.2) -- (0.3,-0.2) node[pos=.5, anchor=north] {\small$\frac{1}{17}$};

\draw (0.3,-0.05) -- (0.3,-0.25);
\draw (0.6,-0.05) -- (0.6,-0.25);
\draw[measline] (0.3,-0.2) -- (0.6,-0.2) node[pos=.5, anchor=north] {\small$\frac{1}{17}$};

\draw (1.1,-0.05) -- (1.1,-0.25);
\draw (1.4,-0.05) -- (1.4,-0.25);
\draw[measline] (1.1,-0.2)-- (1.4,-0.2) node[pos=.5, anchor=north] {\small$\frac{1}{17}$};

\draw (1.4,-0.05) -- (1.4,-0.25);
\draw (2,-0.05) -- (2,-0.25);
\draw[measline] (1.4,-0.2)-- (2,-0.2) node[pos=.5, anchor=north] {\small$\frac{2}{17}$};

\draw [plant] (0.3,0) -- (0.3, 0.3);
\draw [plant] (0.6,0) -- (0.6, 0.3) node[pos=.5, anchor=west] {\small{...}};
\draw [plant] (1.1,0) -- (1.1, 0.3);
\draw [plant] (1.4,0) -- (1.4, 0.3);

\draw[axisline] (-0.2, -0.8 ) -- (0.2, -0.8) node[below] {\small$z$};
\draw[axisline] (-0.2, -0.8) -- (-0.2, -0.4) node[left] {\small$y$};

\end{tikzpicture}

}
	\end{subfigure}
	\caption{Sketch of the configuration for the submerged vegetation example}
	\label{fig:plantpatch_sketch}
\end{figure}
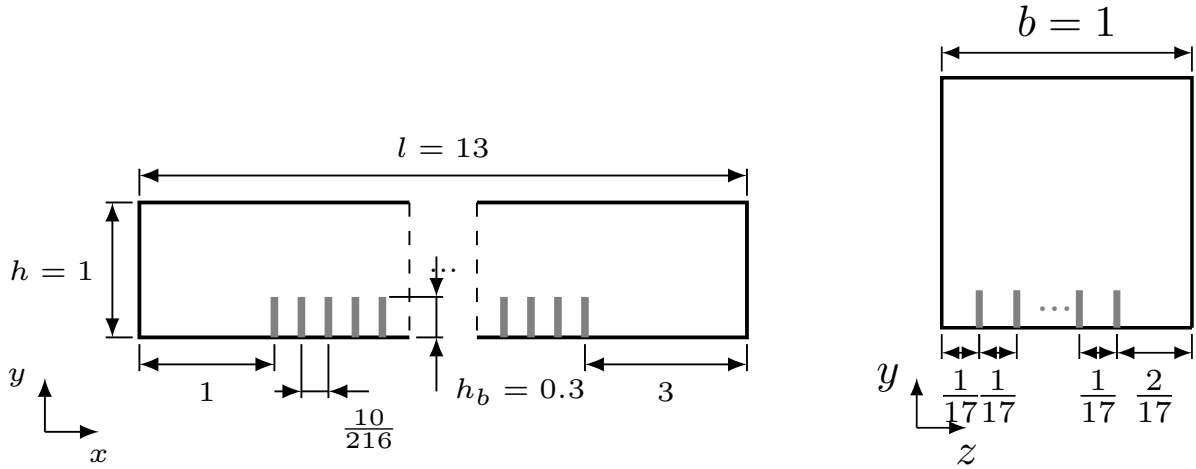

Finally, this example serves as an application of the proposed approach to a problem inspired by a real-life research question: the behavior of a submerged vegetation patch and its impact on the coastal flow.
Flexible terrestrial and aquatic plant canopies massively affect the flow formation in coastal water regions, and in open fields, as their interaction with the flow around them may lead to shear-layer formation and secondary currents. Additionally, plant canopies in water channels and river banks reduce the velocity and throughput of fluid, and can represent a technical challenge concerning the maintenance of the water supply.
Studying and understanding the effect of such submersed vegetation patches is challenging. In~\cite{Okamoto2010, murphy2007}, these phenomena have been studied experimentally and, in~\cite{carollo2005, Ippolito2021, wu1999}, efforts have been made to quantify their effect on the flow around them by analytical solutions to roughness and density variation models.
Recent numerical studies of submerged canopies include the analysis of rigid fixed submerged fibers~\cite{stoesser2009, tariq2022}, 2D models of fluid flow interacting with submerged rod-like  structures, as in~\cite{oconnor2019, oconnor2022, maza2013}, as well as immersed 2D flaps interacting with the 3D flow around it, cf.~\cite{Tschisgale2020, tschisgale2021}. Due to the large number of fibers required to model such canopies realistically, the resolution of such setups in 3D poses a challenge regarding the model complexity. Additionally, the two geometrical scales pose an additional challenge concerning the complexity of the problem setup. Commonly, the displacements of the single beams, which in turn contribute to the patch's overall behavior, are significantly smaller than the channel length. Nevertheless, the behavior of flexible immersed canopies leads to highly interesting flow patterns, making it a prime target application for the proposed FBI approach.

As depicted in Figure~\ref{fig:plantpatch_sketch}, the submerged plant patch is modeled by $225\times15$ slender beams immersed within a three-dimensional fluid channel with the dimensions $13\times1\times1$.
The beams have the length~$h_\indBeam = 0.3$, the radius~$r_\indBeam=10^{-2}$, and are modeled using a hyperelastic material with the Young's modulus~$E^\indBeam=10^7$ and the density~$\rho^\indBeam=10^1$. The beams are fixed to the floor of the channel with equal spacing of~$\dfrac{10}{226}$ in the channel direction, starting at a distance of~$1$ to the inflow boundary. In the direction of the channel's depth, the beams are also spaced equally at an interval of~$\dfrac{1}{17}$, with a distance of~$\dfrac{1}{17}$ to the left wall and~$\dfrac{2}{17}$ to the right wall. This setup leads to a slightly asymmetric behavior, which facilitates flow orthogonal to the channel's principal flow direction.

The fluid is modeled with the density $\rho^\indFluid=1$, the dynamic viscosity $\nu^\indFluid=0.004$, and a no-slip boundary condition is applied to the bottom surface. Non-penetration conditions are applied to the channel surfaces, and a backflow boundary condition is applied to the outflow~\cite{bertoglio2018}. On the inflow boundary, the velocity
\begin{equation*}
	\mathbf{v}_{\mathrm{in}}=100\cdot y\cdot\left(2-y\right),
\end{equation*}
is prescribed in channel direction, and the inflow is zero in both other directions.

The fully-coupled FBI method with mortar-type coupling discretization is applied using linear Lagrange shape functions, a penalty parameter $\epsilon=10^2$, and a time step size $\Delta t=2\cdot 10^{-4}$. The fluid field is discretized with 491,520 finite elements yielding 2,095,236 DoFs, and the 3,375 fibers are discretized with 10,125 finite elements resulting in 81,000 DoFs.

\begin{figure}
	\begin{subfigure}{\textwidth}
		\includegraphics[width=\textwidth]{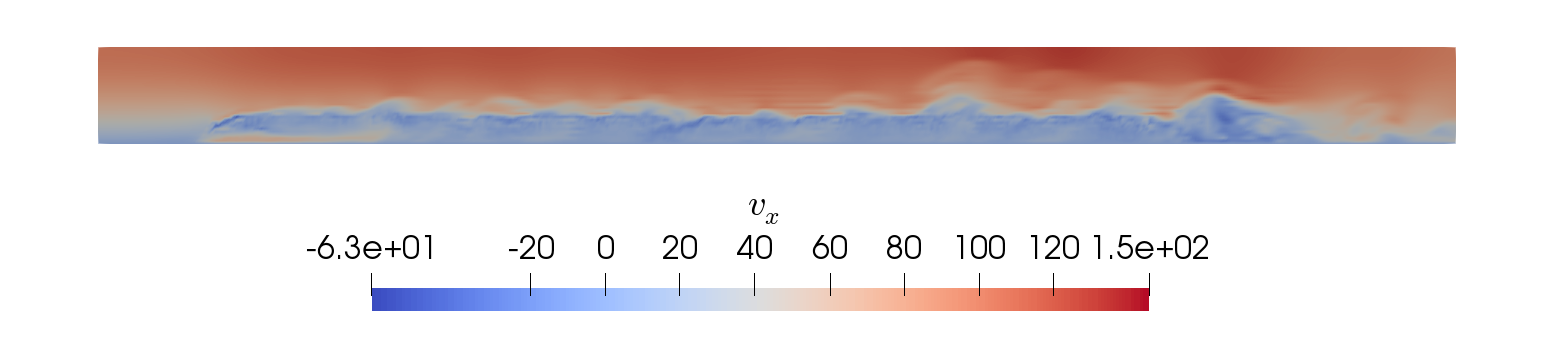}
		\subcaption{Fluid velocity in channel direction on a slice through the channel's middle}
		\label{fig:plantpatch_fluid_xvel}
	\end{subfigure}
	\begin{subfigure}{\textwidth}
		\includegraphics[width=\textwidth]{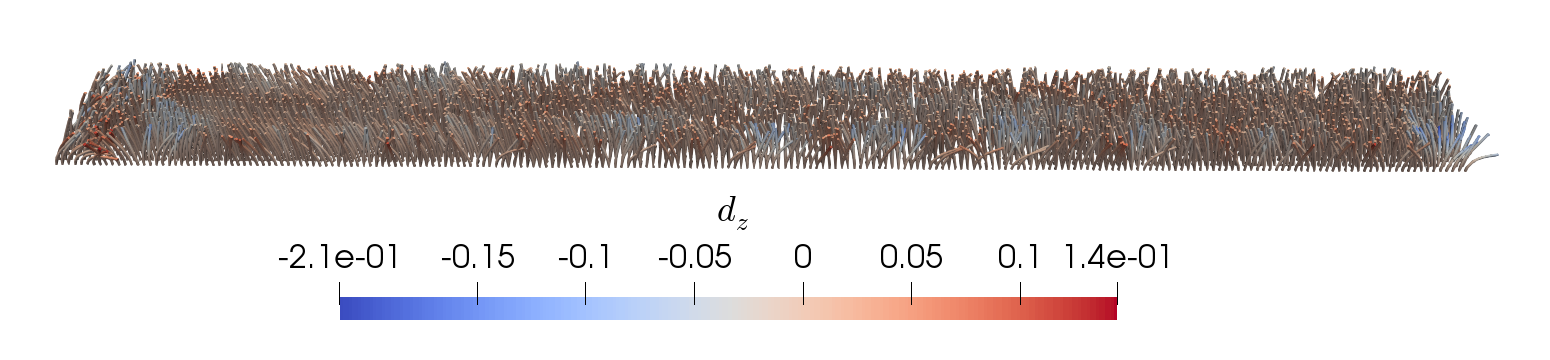}
		\subcaption{Beam displacement orthogonal to the channel length}
		\label{fig:plantpatch_beam_zdisp}
	\end{subfigure}
	\begin{subfigure}{\textwidth}
		\includegraphics[width=\textwidth]{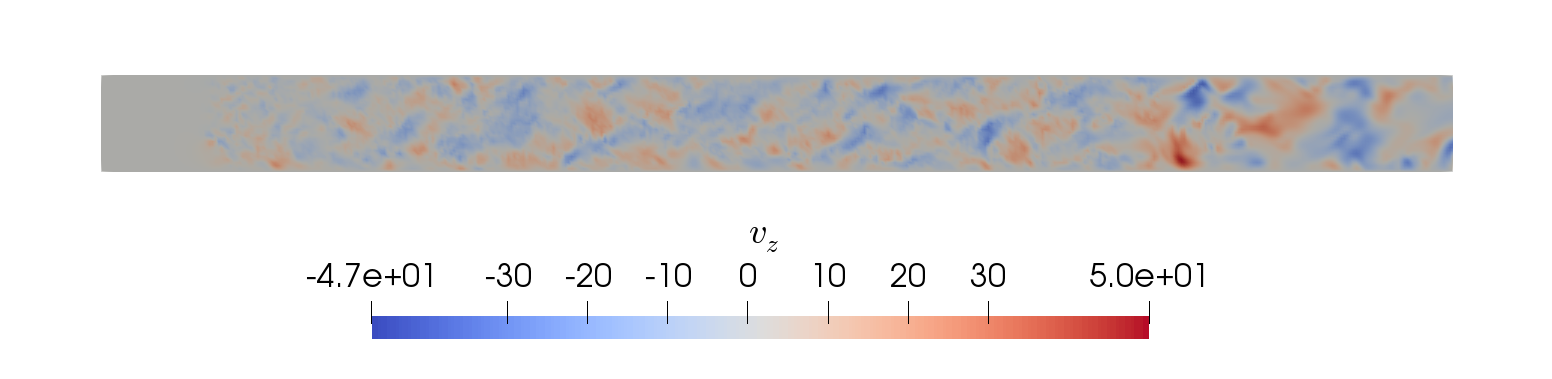}
		\subcaption{Fluid velocity orthogonal to the principal channel direction at height $y=0.3$}
		\label{fig:plantpatch_fluid_zvel}
	\end{subfigure}
	\caption{Plantpatch result}
	\label{fig:plantpatch_results}
\end{figure}

\begin{figure}
	\begin{subfigure}[b]{0.48\textwidth}
		\includegraphics[width=\textwidth]{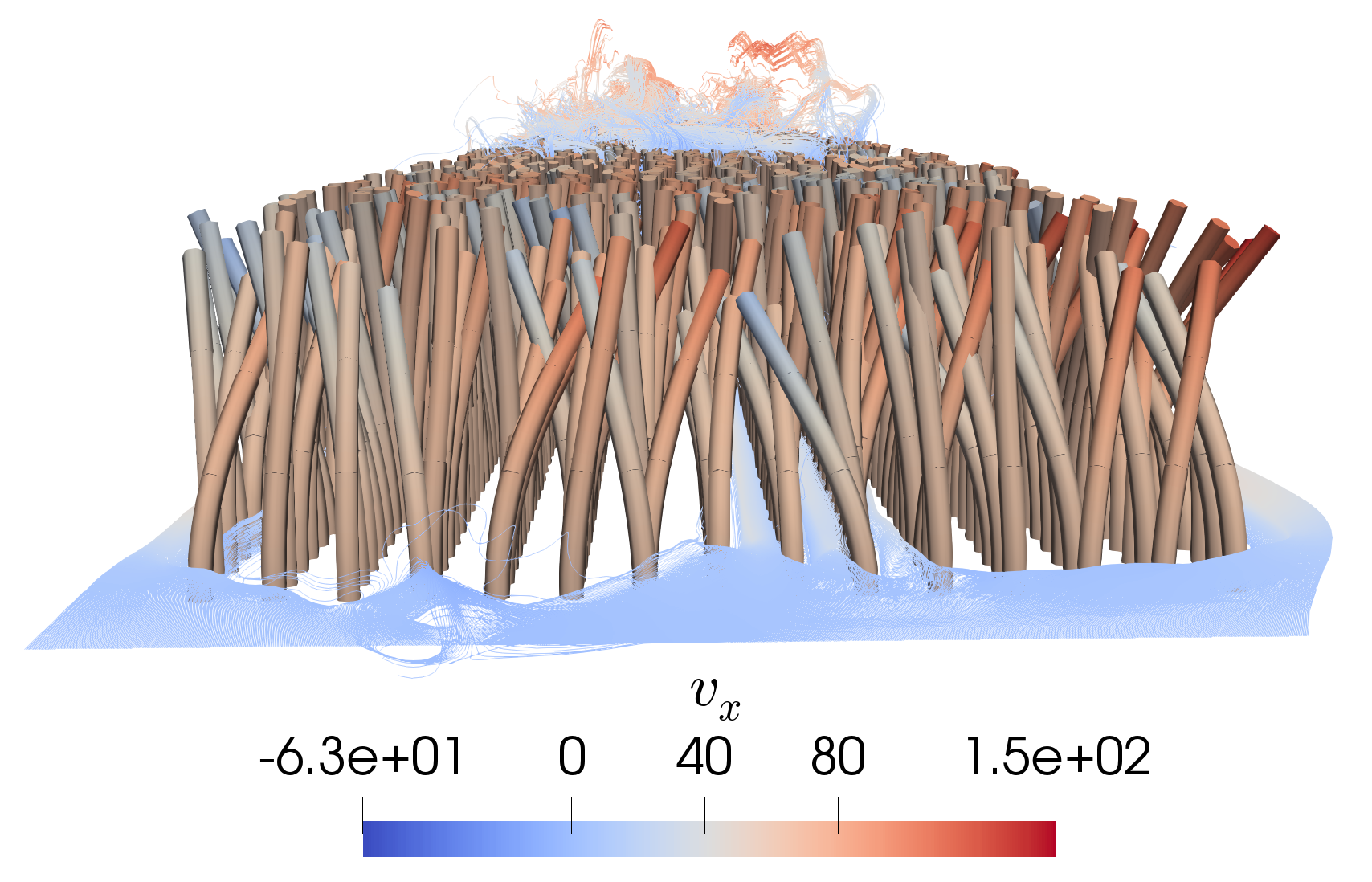}
		\caption{Phenomenon of ground flow}
		\label{fig:plantpatch_zdisp_front_zoom}
	\end{subfigure}
	\begin{subfigure}[b]{0.48\textwidth}
		\includegraphics[width=\textwidth]{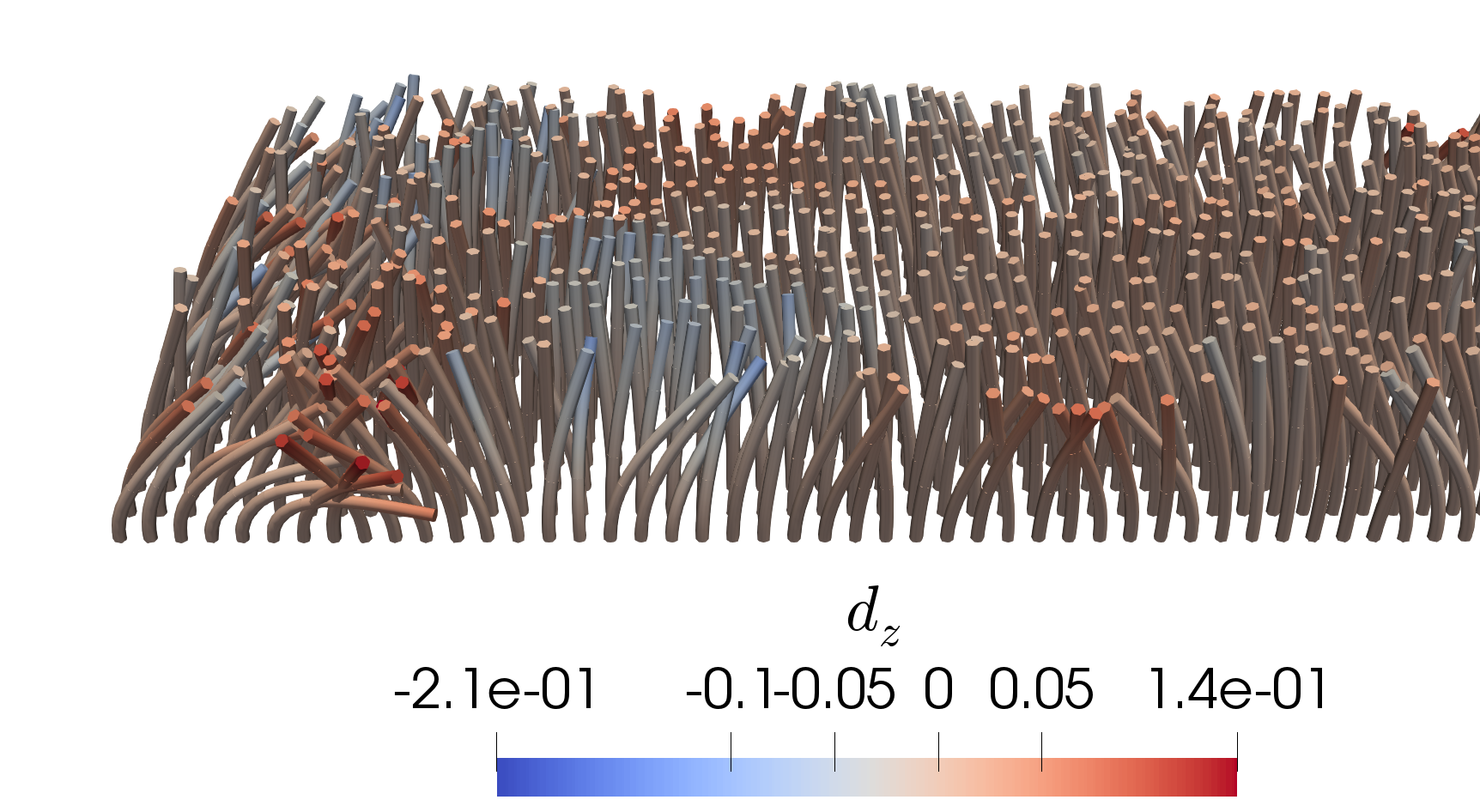}
		\caption{Waves orthogonal to principal flow direction}
		\label{fig:plantpatch_zdisp_side_zoom}
	\end{subfigure}
	\caption{Captured 3D phenomena}
\end{figure}

Figure~\ref{fig:plantpatch_results} depicts the solution of the fluid and beam fields at time $t=0.24$.
The effect of the beams' interaction on the fluid flow is illustrated in Figure~\ref{fig:plantpatch_fluid_xvel}.
The fluid is not only slowed down in the proximity of the fiber patch, but the flexible movement of the beams leads to the formation of monami-type fluid flow as also observed in~\cite{tschisgale2021, Okamoto2010, oconnor2019}.
As depicted in Figure~\ref{fig:plantpatch_beam_zdisp}, this phenomenon goes hand in hand with displacement waves traveling through the beam patch.
This behavior stems from the asymmetry of the model and the successive formation of fluid waves traveling orthogonal to the channel direction, as visualized for the beam displacement $\mathbf{d}$ in Figure~\ref{fig:plantpatch_fluid_zvel}, and also observed in~\cite{tschisgale2021} for immersed flexible flaps.
In contrast to 2D simulations, as reported in~\cite{oconnor2019, oconnor2022}, phenomena caused by flow orthogonal to the principle channel direction can be observed: the fluid pushes the beams in their wake to the sides, effectively increasing the unobstructed flow area in these regions, which in turn leads to variations in the fluid velocity. In contrast to the upper canopy layer of freely moving beam tips, this behavior does not occur near the ground where the beams are fully fixed to the bottom surface. On the contrary, the fixation of the beam near the ground allows the formation of ground flow as hinted at in Figure~\ref{fig:plantpatch_fluid_xvel} and Figure~\ref{fig:plantpatch_zdisp_front_zoom}.

The results of the shown example were obtained without the incorporation of beam-to-beam contact. Although an important phenomenon, note that the beams in the considered example are only subjected to forces applied through interaction with the fluid flow. In this case, when using sufficiently large penalty parameters to enforce the FBI conditions, the beams interact implicitly through the fluid flow. Consequently, the impression of direct beam-to-beam interaction is created, as demonstrated in Figure~\ref{fig:plantpatch_zdisp_side_zoom}. However, a quantitative comparison with the solution of a beam-to-beam scheme such as the ones presented in~\cite{meier_contact, Meier2017} is in order to assess the validity of the presented submerged vegetation model.

This example constitutes the first time that results of a 3D flexible submersed vegetation patch containing structures with high slenderness with respect to two directions are reported. By placing more than 3,000 fibers into the channel, this example showcases an envisioned target application of the proposed FBI approach and demonstrates its potential compared to alternatives, such as fully resolved or homogenized models, when it comes to the modeling of large amounts of fibers and their consequent interaction phenomena. This has been demonstrated by placing a large number of more than $3,000$ stems into the channel.

\section{Conclusions}
\label{sec:conclusions}

We have proposed a computational framework allowing to couple beams with three-dimensional fluid flow via an efficient mixed-dimensional 1D-3D coupling approach. The mesh tying of the beam to the fluid background mesh has been found to yield  good results for mortar-type segment-to-segment coupling using a penalty-regularized Lagrange multiplier method. To solve the resulting surrogate model, a mixed-dimensional multi-physics problem, a matrix-free Newton-Krylov method based on a partitioned Dirichlet-Neumann-type coupling was introduced.

The spatial convergence behavior of the strongly coupled fluid-beam solution under uniform mesh refinement with respect to a fully resolved ALE based 3D FSI solution~\cite{Kloppel2012, mayr2014}
was shown and the approach, thus, quantitatively validated. Further effects of the proposed computational building blocks, such as the mortar-type mesh tying strategy and the acceleration technique for the partitioned solver have been studied with numerical examples.

To demonstrate the robustness and suitability of the proposed approach also for larger examples and application-oriented problems, a fluid channel containing more than $3,000$ submerged fibers, was set up to model the behavior of a submerged vegetation patch. Differences to the solution obtained with existing homogenized, 1D-2D, and 2D-3D coupling strategies have been pointed out.

Further research may relate to the extension of the proposed algorithm to handle combined FBI and classical surface-based FSI interfaces within the same simulation as well as submerged beam-to-beam contact. This will enable the modeling of important effects as they may occur in flow control as well as biomedical problems such as stenting procedures. In particular, the extension to combined FBI and FSI strategies will lead to interesting algorithmic challenges such as the question of the composition of the overall stopping criterion as well as the application of suitable acceleration techniques for both interface types.

\appendix

\section{Exemplary evaluation of mortar matrices}
\label{app:ExampleMortarMatrices}

In the following, the computation of the mortar coupling matrices for a very general embedding of a beam centerline into a linear hexahedral finite element is demonstrated. This example may serve as the basis for test-driven development of a mortar finite element framework for fluid-beam interaction. Here, the positional and tangential degrees of freedom of the beam centerline are given as:
\begin{equation*}
\hat{\mathbf{d}}_h^1 =
\left(\begin{array}{c}
	0.15 \\
	0.2 \\
	0.3
\end{array}
\right), \quad
\hat{\mathbf{d}}_h^2 =
\left(\begin{array}{c}
	0.65 \\
	0.1 \\
	0.1
\end{array}\right), \quad
\hat{\mathbf{t}}_h^1 =
\left(\begin{array}{c}
	0.58 \\
	0.58 \\
	0.58
\end{array}\right), \quad
\hat{\mathbf{t}}_h^2 =
\left(\begin{array}{c}
	0.80 \\
	-0.53 \\
	0.26
\end{array}\right)
\end{equation*}

Furthermore, the positions of the fluid finite element nodes are given as:
\begin{equation*}
\begin{split}
\hat{\mathbf{x}}^{\indFluid,1}_h :=
\left(\begin{array}{c}
-0.95 \\
-0.97 \\
-1.00
\end{array}
\right), \quad
	\hat{\mathbf{x}}^{\indFluid,2}_h :=
	\left(\begin{array}{c}
		0.92 \\
		-1.01 \\
		-1.01
	\end{array}
	\right), \quad
	\hat{\mathbf{x}}^{\indFluid,3}_h :=
	\left(\begin{array}{c}
		0.9 \\
		1.06 \\
		-0.94
	\end{array}
	\right), \quad
	\hat{\mathbf{x}}^{\indFluid,4}_h :=
	\left(\begin{array}{c}
		-1.05 \\
		1.08 \\
		-1.03
	\end{array}
	\right),
\\
	\hat{\mathbf{x}}^{\indFluid,5}_h :=
	\left(\begin{array}{c}
		-1.09 \\
		-1.06 \\
		1.08
	\end{array}
	\right)
, \quad
	\hat{\mathbf{x}}^{\indFluid,6}_h :=
	\left(\begin{array}{c}
		0.97 \\
		-1.01 \\
		0.92
	\end{array}
	\right)
, \quad
	\hat{\mathbf{x}}^{\indFluid,7}_h :=
	\left(\begin{array}{c}
		1.09 \\
		1.03 \\
		0.96
	\end{array}
	\right)
, \quad
	\hat{\mathbf{x}}^{\indFluid,8}_h :=
	\left(\begin{array}{c}
		-0.94 \\
		0.95 \\
		0.96
	\end{array}
	\right)
	\end{split}
\end{equation*}

Linear Lagrange shape functions are used for the 8-noded hexahedral element:
\begin{equation*}
	\begin{split}
	N_1 = \dfrac{1}{8}\left(1-\xi_1\right)\left(1-\xi_2\right)\left(1-\xi_3\right), \quad
	N_2 = \dfrac{1}{8}\left(1+\xi_1\right)\left(1-\xi_2\right)\left(1-\xi_3\right), \\
	N_3 = \dfrac{1}{8}\left(1+\xi_1\right)\left(1+\xi_2\right)\left(1-\xi_3\right), \quad
	N_4 = \dfrac{1}{8}\left(1-\xi_1\right)\left(1+\xi_2\right)\left(1-\xi_3\right), \\
	N_5 = \dfrac{1}{8}\left(1-\xi_1\right)\left(1-\xi_2\right)\left(1+\xi_3\right), \quad
	N_6 = \dfrac{1}{8}\left(1+\xi_1\right)\left(1-\xi_2\right)\left(1+\xi_3\right), \\
	N_7 = \dfrac{1}{8}\left(1+\xi_1\right)\left(1+\xi_2\right)\left(1+\xi_3\right), \quad
	N_8 = \dfrac{1}{8}\left(1-\xi_1\right)\left(1+\xi_2\right)\left(1+\xi_3\right).
	\end{split}
\end{equation*}

The 3rd order Hermite shape functions suitable for $C^1$-continuous beam centerlines\cite{Meier2016_contact} are used for the spatial discretization of the beam centerline, reading
\begin{equation*}
	\begin{split}
H_d^1 =  \dfrac{1}{4}\left(2+\xi\right)\left(1-\xi\right)^2, \quad H_d^2 = \dfrac{1}{4}\left(2-\xi\right)\left(1+\xi\right)^2, \\
H_t^1 =  \dfrac{1}{4}\left(1+\xi\right)\left(1-\xi\right)^2, \quad H_t^2 = \dfrac{1}{4}\left(1-\xi\right)\left(1+\xi\right)^2,
\end{split}
\end{equation*}

while the linear Lagrange shape functions for the discretization of the Lagrange multiplier are given by
\begin{equation*}
\Phi_1 =  \dfrac{1}{2}\left(1+\xi\right), \quad \Phi_2 =  \dfrac{1}{2}\left(1-\xi\right).
\end{equation*}

In order to compute the mortar coupling matrices~\eqref{eq:mortarmatrices} via numerical integration,
the projection of all beam centerline Gauss points onto the fluid element according to Remark~\ref{remark:segmentation} have to be found.
More information on the projection and segmentation procedure can be found in~\cite{steinbrecher2020}.

For the sake of readability, we will drop rows and columns of the fluid-sided matrix~$\mathbf{D}$ that are related to pressure degrees of freedom.
Since they are excluded from the mesh coupling by design, they would carry zeroes only anyways.
This reduces the matrix size from $6\times32$ to $6\times12$ in this example.
On that basis, the mortar coupling matrices, where the values are given for four positions after the decimal point, can be computed as:
\begin{equation*}
	\mathbf{D} =
	\left(
	\begin{array}{cccccccccccc}
		0.1954 & 0 & 0  & 0.01819 & 0 & 0 & 0.0989 & 0 & 0 & -0.0135 & 0 & 0\\
		0 & 0.1954 & 0 & 0 & 0.01819 & 0 & 0 & 0.0989 & 0 & 0 & -0.0135 & 0 \\
		0 & 0 & 0.1954 & 0 & 0 & 0.01819 & 0 & 0 & 0.0989 & 0 & 0 & -0.0135 \\
		0.0947 & 0 & 0 & 0.0135 & 0 & 0 & 0.2301 & 0 & 0 & -0.0208 & 0 & 0 \\
		0 & 0.0947 & 0 & 0 & 0.0135 & 0 & 0 & 0.2301 & 0 & 0 & -0.0208 & 0 \\
		0 & 0 & 0.0947 & 0 & 0 & 0.0135 & 0 & 0 & 0.2301 & 0 & 0 & -0.0208
	\end{array}
    \right)
\end{equation*}

\begin{equation*}
	\begin{split}
	\mathbf{M} =
	&\left(
	\begin{array}{cccccccccccc}
	0.0140 & 0 & 0 & 0.0282 & 0 & 0 & 0.0433 & 0 & 0 & 0.0218 & 0 & 0 \\
	0 & 0.0140 & 0 & 0 & 0.0282 & 0 & 0 & 0.0433 & 0 & 0 & 0.0281 & 0 \\
	0 & 0 & 0.0140 & 0 & 0 & 0.0282 & 0 & 0 & 0.0433 & 0 & 0 & 0.0218 \\
	0.0137 & 0 & 0 & 0.0417 & 0 & 0 & 0.0581 & 0 & 0 & 0.0198 & 0 & 0 \\
	0 & 0.0137 & 0 & 0 & 0.0417 & 0 & 0 & 0.0581 & 0 & 0 & 0.0198 & 0 \\
	0 & 0 & 0.0137 & 0 & 0 & 0.0417 & 0 & 0 & 0.0581 & 0 & 0 & 0.0198
    \end{array}
	\right.
	\\
	\\
	&\left.
		\begin{array}{cccccccccccc}
		0.0250 & 0 & 0 & 0.0482 & 0 & 0 & 0.0747 & 0 & 0 & 0.0391 & 0 & 0 \\
		0 & 0.0250 & 0 & 0 & 0.0482 & 0 & 0 & 0.0747 & 0 & 0 & 0.0391 \\
		0 & 0 & 0.0250 & 0 & 0 & 0.04822 & 0 & 0 & 0.0747 & 0 & 0 & 0.0391 \\
		0.0203 & 0 & 0 & 0.0587 & 0 & 0 & 0.0829 & 0 & 0 & 0.0296 & 0 \\
		0 & 0.0203 & 0 & 0 & 0.0587 & 0 & 0 & 0.0829 & 0 & 0 & 0.0296 & 0 \\
		0 & 0 & 0.0203 & 0 & 0 & 0.0587 & 0 & 0 & 0.0829 & 0 & 0 & 0.0296
	\end{array}
	\right)
	\end{split}
\end{equation*}

\bibliographystyle{abbrv}
\bibliography{references}%

\end{document}